\documentclass[preprint,12pt]{elsarticle}
\usepackage{amssymb}
\usepackage{bbm}
\usepackage{float}
\usepackage{natbib}
\usepackage{amsmath}
\usepackage{nomencl}
\makenomenclature
\usepackage{caption}
\usepackage{subcaption}
\usepackage[export]{adjustbox}
\usepackage[colorlinks,linkcolor=red,anchorcolor=blue,citecolor=green]{hyperref}
\usepackage{geometry}
\usepackage{mathrsfs}
\usepackage{cleveref}

\usepackage{amsthm}

\journal{Elsevier}

\begin{document}

\begin{frontmatter}

\title{A Gas-Kinetic Scheme for Maxwell Equations}

\author[a]{Zhigang PU}
\ead{zgpuac@ust.hk}
\author[a]{Wenpei Long}
\ead{wlongab@connect.ust.hk}
\author[a,b]{Kun XU\corref{cor1}}
\ead{makxu@ust.hk}
\cortext[cor1]{Cooresponding author}

\address[a]{{Department of Mathematics, Hong Kong University of Science and Technology},
            {Clear Water Bay, Kowloon},
            {Hong Kong},
            {China}}
\address[b]{{Shenzhen Research Institute, Hong Kong University of Science and Technology},
            {Shenzhen},
            {China}}
\begin{abstract}
The Gas-Kinetic Scheme (GKS), widely used in computational fluid dynamics for simulating hypersonic and other complicated flow phenomena, is extended in this work to electromagnetic problems by solving Maxwell's equations. In contrast to the classical GKS formulation, the proposed scheme employs a discrete rather than a continuous velocity space. By evaluating a time-accurate numerical flux at cell interfaces, the proposed scheme attains second-order accuracy within a single step. Its kinetic formulation provides an inherently multidimensional framework, while the finite-volume formulation ensures straightforward extension to unstructured meshes. Through the incorporation of a collision process, the scheme exhibits lower numerical dissipation than classical flux-vector splitting (FVS) methods.
Furthermore, the kinetic decomposition enables direct implementation of non-reflecting boundary conditions. The proposed scheme is validated against several benchmark problems and compared with established methods, including the Finite-Difference Time-Domain (FDTD) method and FVS. A lattice Boltzmann method (LBM) implementation is also included for comparative analysis. Finally, the technique is applied to simulate electromagnetic wave propagation in a realistic aircraft configuration, demonstrating its ability to model complex geometries.
\end{abstract}


\end{frontmatter}

\printnomenclature

\section{Introduction}
\label{introduction}

Numerical methods for solving Maxwell's equations are fundamental to computational electromagnetics (CEM). The development of accurate and efficient numerical methods serves both practical applications and fundamental research needs. These methods are crucial for the design and modeling of communication systems, including antennas, radars, and satellites. The most commonly used approaches are the finite-difference time-domain (FDTD) method, finite element method (FEM), and finite volume method (FVM), each with distinct advantages and limitations.

The FDTD method is widely utilized in computational electromagnetics for its simplicity and efficiency. It employs staggered Yee grid arrangements to satisfy elliptic divergence constraints, with electric and magnetic fields defined at staggered grid points and updated using a leapfrog scheme. This approach achieves second-order accuracy in both time and space \cite{taflove2005computational,yee1966numerical}. While Holland extended the Yee algorithm to curvilinear structured mesh \cite{holland1983finite}, the FDTD method faces challenges when applied to unstructured meshes for complex geometries \cite{yee_finite-difference_1997}.
FEM offers greater flexibility in handling complex geometries than FDTD. Based on weak formulations of Maxwell equations in the frequency domain, FEM uses edge elements introduced by Nédélec to address spurious modes caused by nodal basis functions \cite{nedelec1980mixed}. These elements allow discontinuities in normal components across interfaces while maintaining tangential component continuity. A limitation of FEM is its requirement to solve a sparse linear system at each time step, though techniques like mass lumping have been developed to simplify the mass matrix to diagonal form \cite{fisher2005generalized}. Similar to FEM, the method of moments (MoM) solves the integral form of Maxwell equations \cite{harrington1982field,bondeson2012computational}. For materials with spatially varying electromagnetic properties or inhomogeneities, differential Maxwell systems are preferred over integral equations.

FVM, initially developed for nonlinear hyperbolic systems in fluid dynamics \cite{jameson1981numerical,leveque2002finite}, integrates equations over mesh cells where cell-averaged quantities evolve based on interface fluxes calculated through Riemann problems \cite{toro2013riemann}. Mohammadian achieved second-order accuracy for Maxwell equations using the Lax-Wendroff scheme \cite{mohammadian1991computation}, while Shang obtained higher-order spatial accuracy using flux vector splitting and the MUSCL scheme \cite{shang_characteristic-based_1995}. The characteristic-based Riemann solution effectively handles boundary conditions for infinitely outgoing waves, preventing spurious reflections. Munz introduced perfect hyperbolic Maxwell (PHM) equations to control divergence error in particle-in-cell (PIC) simulations and later integrated PHM into the FVM framework \cite{munz_divergence_2000}.
The discontinuous Galerkin (DG) method combines FVM and FEM to avoid complex reconstruction and to achieve higher spatial accuracy \cite{taube_high-order_2009}. While DG with explicit Runge-Kutta time integration is subject to restrictive CFL conditions \cite{toulorge2011cfl}, implicit methods, such as the kinetic implicit scheme, help overcome these limitations \cite{gerhard_unconditionally_2022}. In computational electromagnetics, the finite-volume formulation of the integral form of Maxwell's equations is known as the finite integration method (FIT) \cite{weiland1977discretization}. Initially developed independently of FDTD and later expanded as its generalization \cite{weiland1984numerical,weiland1996time}, FIT employs staggered grids but uses numerical flux based on the integral form of Maxwell's equation rather than Riemann solutions. FIT provides greater flexibility in grid selection for complex geometries and eliminates the need for reconstruction through dual grids \cite{carstensen2012computational}.

Since both Maxwell's equations and fluid dynamics involve hyperbolic systems, methods from computational fluid dynamics (CFD) can be adapted for electromagnetic computations. Beyond Riemann-solver-based FVM methods, kinetic-solver-based approaches provide concise formulations for multidimensional hyperbolic systems without requiring local Riemann problem solutions \cite{aregba2000discrete}. In fluid dynamics, prominent kinetic-theory-based methods include the Gas-Kinetic Scheme (GKS) \cite{xu2001gas}, which uses a continuous velocity space, and the Lattice Boltzmann Method (LBM) \cite{guo2013lattice}, which employs a discrete velocity space. These kinetic models offer simpler formulations for multidimensional systems compared to Riemann-solver-based schemes.
In GKS, the integral solution of the Bhatnagar–Gross–Krook (BGK) equation at cell interfaces extends directly to multiple dimensions, while LBM achieves true multidimensionality through an isotropic lattice velocity distribution and precise transport operators. Recently, there have been several works applying the vector lattice Boltzmann method (LBM) to Maxwell equations \cite{mendoza2010three, hanasoge2011lattice, liu2014lattice,hauser2019comparison}. The multidimensional nature of lattice schemes inherently enhances the preservation of elliptical divergence constraints, as demonstrated by Dellar et al. \cite{dellar_lattice_2002}. Khan et al. \cite{khan2024electromagnetic} explored the efficacy of the LBM in simulating electromagnetic wave scattering from curved and complex surfaces. Gerhard's research approximates Maxwell's equations using discrete kinetic models solved via Discontinuous Galerkin (DG) methods \cite{gerhard_unconditionally_2022,gerhard2024parallel}. This approach facilitates the use of sweep algorithms to achieve large time steps, thereby circumventing the need to solve large linear systems inherent in classical implicit schemes.

This work extends the Gas-Kinetic Scheme (GKS) to solve Maxwell's equations, representing its first application to electromagnetic problems. The GKS framework has seen considerable success in compressible fluid dynamics. It recovers the Navier-Stokes equations in the continuum limit and supports later milestone methods like the Unified Gas-Kinetic Scheme (UGKS) and the Unified Gas-Kinetic Wave-Particle (UGKWP) method for rarefied flows and plasma simulations \cite{xu2010unified, liu2021unified, pu_gas-kinetic_2024, pu2025unified}. By adapting GKS to electromagnetics, we are able to establish a consistent framework for plasma simulation. The proposed scheme uses a time-accurate flux evaluation at cell interfaces. This allows it to achieve second-order accuracy in a single time step. Future work will exploit this feature to develop compact high-order solvers for Maxwell's equations \cite{zhao2019compact, ji2020hweno, ji2021compact}. The finite-volume formulation permits easy extension to unstructured meshes.
Furthermore, the kinetic basis of the scheme provides a naturally multidimensional numerical flux. We validate the scheme against standard benchmark problems. Performance is compared to established methods, including the FDTD method and FVS solver. A LBM solver is also included for comparison. Finally, we demonstrate the method's capability for complex geometries by simulating electromagnetic wave propagation around a realistic aircraft configuration.

This paper is organized as follows. Section \ref{sec:kinetic representation} introduces the kinetic representation of the Maxwell system. Section \ref{sec:algorithms} presents the numerical algorithms, including the proposed GKS, as well as the LBM and FVS methods used for comparison. Section \ref{sec:analysis} provides an analysis of the proposed schemes. Numerical studies and test cases are presented in Section \ref{sec:numerical studies}. Finally, Section \ref{sec:conclusion} concludes the paper.

\section{Kinetic representation of Maxwell system}
\label{sec:kinetic representation}

\subsection{Maxwell system}
The PHM equations are used here for divergence cleaning,
\begin{align}
    & \frac{\partial \boldsymbol{E}}{\partial t}-c^2 \nabla_{\boldsymbol{x}} \times \boldsymbol{B}+\chi c^2 \nabla_{\boldsymbol{x}} \phi=-\frac{1}{\epsilon_0} \boldsymbol{J},\label{eq:PHM Amphere law}\\
    & \frac{\partial \boldsymbol{B}}{\partial t}+\nabla_{\boldsymbol{x}} \times \boldsymbol{E}+\gamma \nabla_{\boldsymbol{x}} \psi=0, \label{eq:PHM Faraday law}\\
    & \frac{1}{\chi} \frac{\partial \phi}{\partial t}+\nabla_{\boldsymbol{x}} \cdot \boldsymbol{E}=\frac{\rho}{{\epsilon_0}}, \label{eq:PHM E divergence}\\
    & \frac{\epsilon_0 \mu_0}{\gamma} \frac{\partial \psi}{\partial t}+\nabla_{\boldsymbol{x}} \cdot \boldsymbol{B}=0,\label{eq:PHM B divergence}
\end{align}
where $\phi,\psi$ are artificial correction potentials to accommodate divergence errors traveling at speed $\gamma c$ and $\chi c$ \cite{munz_divergence_2000}. To construct a finite-volume scheme, we rewrite Eqs.\eqref{eq:PHM Amphere law}, \eqref{eq:PHM Faraday law}, \eqref{eq:PHM E divergence} and \eqref{eq:PHM B divergence} as a system of linear hyperbolic evolution equations
\begin{equation}
    \frac{\partial \boldsymbol{U}}{\partial t} + \nabla\cdot\mathbbm{F}(\boldsymbol{U}) = \boldsymbol{S}(\boldsymbol{U}),
    \label{eq:conservation law}
\end{equation}
where  $\boldsymbol{U}\in\mathbb{R}^8$ is the vector of unknowns given by
\begin{equation}
    \boldsymbol{U}=\left(E_{x}, E_{y}, E_{z}, B_{x}, B_{y}, B_{z}, \phi, \psi\right)^{T},
\end{equation}
and $\mathbbm{F}(\boldsymbol{U}) = (\boldsymbol{F}_1(\boldsymbol{U}), \boldsymbol{F}_2(\boldsymbol{U}), \boldsymbol{F}_3(\boldsymbol{U}))\in \mathbb{R}^{8\times 3}$ is the flux tensor. The flux function along $j$ direction $\boldsymbol{F}_j(\boldsymbol{U})=\mathbbm{A}_j\boldsymbol{U}$, the Jacobian matrixs $\mathbbm{A}_j\in\mathbb{R}^{8\times 8}$ with constant entries are defined as
\begin{equation}
\mathbbm{A}_j=\left(\begin{array}{cccccccc}
0 & 0 & 0 & & &  & \chi c^2 \delta_{1 j} &0\\
0 & 0 & 0 & & c^2 \mathbbm{M}_j & & \chi c^2 \delta_{2 j} &0\\
0 & 0 & 0 & & & & \chi c^2 \delta_{3 j} &0\\
& &   & 0 & 0 & 0 & 0 &\gamma \delta_{1 j} \\
& \mathbbm{M}_j^T &   & 0 & 0 & 0 & 0 & \gamma \delta_{2 j} \\
& &   & 0 & 0 & 0 & 0 & \gamma \delta_{3 j}\\
\chi \delta_{1 j} & \chi \delta_{2 j} & \chi \delta_{3 j} & 0 & 0 & 0 & 0 & 0  \\
0 & 0 & 0 &\gamma c^2 \delta_{1 j} & \gamma c^2\delta_{2 j} & \gamma c^2 \delta_{3 j} & 0 & 0
\end{array}\right) ; \quad j=1,2,3,
\label{eq:Jacobian}
\end{equation}
where $\delta_{i j}$ denotes the usual Kronecker symbol and the two $3 \times 3$ matrices $\mathbbm{M}_j$ are found to be
$$
\mathbbm{M}_1=\left(\begin{array}{ccc}
0 & 0 & 0 \\
0 & 0 & 1 \\
0 & -1 & 0
\end{array}\right), \quad \mathbbm{M}_2=\left(\begin{array}{ccc}
0 & 0 & -1 \\
0 & 0 & 0 \\
1 & 0 & 0
\end{array}\right),\quad \mathbbm{M}_3=\left(\begin{array}{ccc}
0 & 1 & 0 \\
-1 & 0 & 0 \\
0 & 0 & 0
\end{array}\right).
$$
Finally, the flux function in different directions is given as
$$
\boldsymbol{F}_1(\boldsymbol{U})
\begin{pmatrix}
\chi c^2 \phi \\
c^2 B_z \\
-c^2 B_y \\
\gamma \psi \\
-E_z \\
E_y \\
\chi E_x \\
\gamma c^2 B_x
\end{pmatrix}\quad \boldsymbol{F}_2(\boldsymbol{U}) = 
\begin{pmatrix}
-c^2 B_z \\
\chi c^2 \phi \\
c^2 B_x \\
E_z \\
\gamma \psi \\
-E_x \\
\chi E_y \\
\gamma c^2 B_y
\end{pmatrix}\quad \boldsymbol{F}_3(\boldsymbol{U}) = 
\begin{pmatrix}
c^2 B_y \\
-c^2 B_x \\
\chi c^2 \phi \\
-E_y \\
E_x \\
\gamma \psi \\
\chi E_z \\
\gamma c^2 B_z
\end{pmatrix}.
$$
The source term reads as
\begin{equation}
 \boldsymbol{S}(\boldsymbol{U})=\left(-J_{x} / \epsilon_0,-J_{y} / \epsilon_0,-J_{z} / \epsilon_0, 0,0,0, \chi \rho / \epsilon_{0}, 0\right)^{T}.
\end{equation}

\subsection{Kinetic representation}

Following Dellar \cite{dellar_lattice_2002} and Bouchut \cite{bouchut_construction_nodate}, we introduce a vector-valued distribution function $\boldsymbol{f}(\boldsymbol{x},\boldsymbol{u}_k,t)\in \mathbb{R}^{8\times 1}$ for the system,
\begin{equation}
    \sum_{k=1}^M \boldsymbol{f}(\boldsymbol{x},\boldsymbol{u}_k,t)  = \boldsymbol{U}(\boldsymbol{x},t),
    \label{eq:def_f_discretized}
\end{equation}
where $\boldsymbol{u}_k = (u_k, u_k, w_k) \in\mathbb{R}^{3}$ is the $k$th discrete kinetic velocity vector, $u_k, u_k, w_k$ is the x, y, z velocity component, $M$ is the total number of discrete kinetic velocities. The discretized BGK kinetic equation is given as:
\begin{equation}
    \frac{\partial \boldsymbol{f}_k}{\partial t} + \boldsymbol{u}_k\cdot\nabla\boldsymbol{f}_k = \frac{\boldsymbol{g}_k-\boldsymbol{f}_k}{\tau},
    \label{eq:discretized vector BGK}
\end{equation}
where $\boldsymbol{g}_k=\boldsymbol{g}(\boldsymbol{x},\boldsymbol{u}_k,t)\in\mathbb{R}^{8}$ is the LTE distribution function corresponding to the discrete velocity $\boldsymbol{u}_k$. To find the expression of $\boldsymbol{g}_k$, in the limit $\tau\rightarrow 0$  the microscopic system Eq.\eqref{eq:discretized vector BGK} should recover the macroscopic system \eqref{eq:conservation law}. In such a limit, with the local equilibrium assumption $\boldsymbol{f}_k=\boldsymbol{g}_k$, the  distribution must satisfy the following compatibility constraints
\begin{equation}
    \sum_{k=1}^M \boldsymbol{g}_k = \boldsymbol{U},\quad
    \sum_{k=1}^M u_k^j \boldsymbol{g}_k = \boldsymbol{F}_j(\boldsymbol{U}),
\label{eq: constraints of LTE}
\end{equation}
where the dimensional index $j=1,2,3$. For different components of $\boldsymbol{u}_k = (u_k^1, u_k^2, u_k^3)^T$, we also denoted as $\boldsymbol{u}_k = (u_k, u_k, w_k)^T$. For simplicity, the summation symbols are simplified as 
$$
\begin{aligned}
<\boldsymbol{g}> = \sum_{k=1}^M \boldsymbol{g}_k, \quad <\boldsymbol{g}>_{>0} = \sum_{\boldsymbol{u}_k\cdot \boldsymbol{n}_s > 0} \boldsymbol{g}_k, \quad <\boldsymbol{g}>_{<0} = \sum_{\boldsymbol{u}_k\cdot \boldsymbol{n}_s < 0} \boldsymbol{g}_k,\quad  <u\boldsymbol{g}> = \sum_{k=1}^M u\boldsymbol{g}_k.
\end{aligned}
$$
In this work, $M$ is equal to $4$ to satisfy the $4$ constraints for $\boldsymbol{g}_k$ in Eq.\eqref{eq: constraints of LTE}.
The four kinetic velocities are chosen as \cite{gerhard_unconditionally_2022}
\begin{equation}
    \boldsymbol{u}_0=\left(\begin{array}{l}1 \\1 \\1\end{array}\right)c, \quad \boldsymbol{u}_1=\left(\begin{array}{c}1 \\-1 \\-1\end{array}\right)c, \quad \boldsymbol{u}_2=\left(\begin{array}{c}-1 \\1 \\-1\end{array}\right)c, \quad \boldsymbol{u}_3=\left(\begin{array}{c}-1 \\-1 \\1\end{array}\right)c,
    \label{eq:kinetic velocities in 3D}
\end{equation}
where $c$ is the speed of light. The LTE coefficients under such conditions can be calculated as
\begin{equation}
    \boldsymbol{g}_k = \frac{\boldsymbol{U}_k}{4} + \frac{1}{4c^2}\mathbbm{F}(\boldsymbol{U})\cdot\boldsymbol{u}_k.
    \label{eq:eq for 3D}
\end{equation}

\begin{figure}
    \centering
    \begin{subfigure}{0.32\textwidth}
        \centering
        \includegraphics[width=1.1\textwidth]{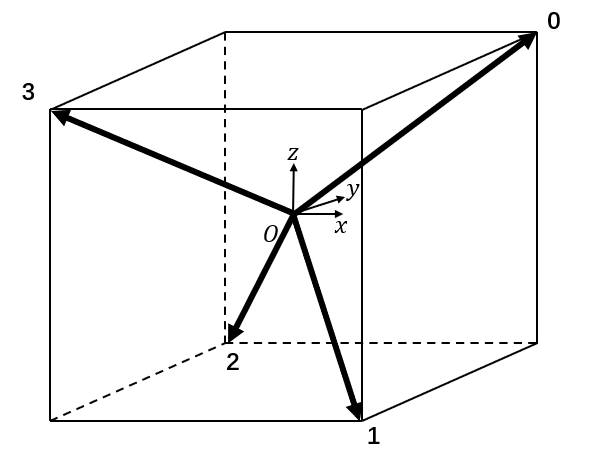}
    \end{subfigure}
    \begin{subfigure}{0.32\textwidth}
        \centering
        \includegraphics[width=0.85\textwidth]{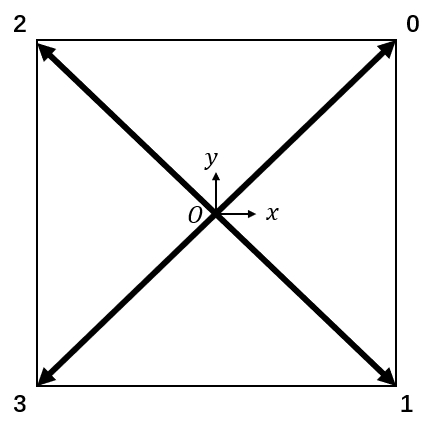}
    \end{subfigure}
    \begin{subfigure}{0.32\textwidth}
        \centering
        \includegraphics[width=0.85\textwidth]{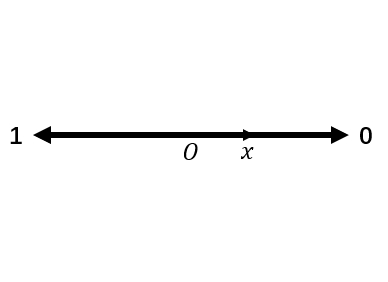}
    \end{subfigure}
    \caption{The lattice structure of D3Q4, D2Q4 and D1Q2.}
    \label{fig:lattice structure}
\end{figure}
This lattice automatically degenerates to D2Q4 in 2D and D1Q2 in 1D as shown in Figure \ref{fig:lattice structure}.

\section{Algorithms}
\label{sec:algorithms}

\subsection{GKS method}
In the framework of the Finite Volume Method (FVM), the cell-averaged variables $\boldsymbol{U}_i$ in a physical cell $\Omega_i$ are defined as 
$$
\boldsymbol{U}_i = \frac{1}{|\Omega_i|} \int_{\Omega_i} \boldsymbol{U}(\boldsymbol{x}) \mathrm{d} \boldsymbol{x},
$$
where $|\Omega_i|$ is the volume of cell $\Omega_i$. For a discretized time step $\Delta t = t^{n+1} - t^n$, the evolution of $\boldsymbol{U}_i$ is 
$$
\boldsymbol{U}_i^{n+1} = \boldsymbol{U}_i^{n} - \frac{\Delta t}{|\Omega_i|} \sum_{s\in\partial \Omega_i} |l_s| \mathscr{F}_{s} + \Delta t \boldsymbol{S}_i^n,
$$
where $l_s$ is the cell interface with center $\boldsymbol{x}_s$ and outer unit normal vector $\boldsymbol{n}_s$. $|l_s|$ is the area of the cell interface. $\boldsymbol{S}_i$ is the source term due to the current density and charge density. The numerical flux across the interface $\mathscr{F}_{s}$ can be evaluated by coordinate transformation,
$$
\mathscr{F}_{s}\cdot\boldsymbol{n}_s =\mathbbm{T}^{-1}[\widetilde{\mathscr{F}}_s(\mathbbm{T}\boldsymbol{U})\cdot\widetilde{\boldsymbol{n}}_s] =  \mathbbm{T}^{-1}[\widetilde{\mathscr{F}}_s(\widetilde{\boldsymbol{U}})\cdot\widetilde{\boldsymbol{n}}_s],
$$
where $\mathbbm{T}=\text{diag}(\mathbbm{T}^{'},\mathbbm{T}^{'},1,1)$ is a rotating operator, transforming the global coordinate system into a local coordinate system, and
\begin{equation*}
\mathbbm{T}^{\prime}=\left(\begin{array}{ccc}
n_1 & n_2 & n_3 \\
-n_2 & n_1+\frac{n_3^2}{1+n_1} & -\frac{n_{2} n_3}{1+n_1} \\
-n_3 & -\frac{n_{2} n_3}{1+n_1} & 1-\frac{n_3^2}{1+n_1}
\end{array}\right), \quad n_1 \neq-1,
\end{equation*}
when $n_1=-1$, $\mathbbm{T}^{'} = \text{diag}(-1,-1,1)$. $\widetilde{\mathscr{F}}_s$ means numerical flux in the local coordinate system, $\widetilde{\boldsymbol{U}}$ means variables in the local coordinate system. Operator $\mathbbm{T}^{-1}$ transforms the local coordinate system to the global coordinate system. For simplicity, The tilde symbol is omitted in later discussion.

Numerical flux can be evaluated from the distribution function at the interface,
$$
\mathscr{F}_{s} = \frac{1}{\Delta t} \int_{t^n}^{t^{n+1}} \sum_{k=1}^M \boldsymbol{u}_k \cdot \boldsymbol{n}_s \boldsymbol{f}_k (\boldsymbol{x}_s, t) dt.
$$
The distribution function on the cell interface can be calculated by the integral solution of the BGK equation \ref{eq:discretized vector BGK}. For simplicity, the subscript $k$ is omitted here. To construct the numerical flux at $\boldsymbol{x}=(0,0,0)^T$, the time-dependent solution on the interface can be written as,
\begin{equation}
\boldsymbol{f}_k(\boldsymbol{x},t) = \frac{1}{\tau} \int_0^t \boldsymbol{g}_{0k}(\boldsymbol{x}^{'},t^{'})e^{-(t-t^{'})/\tau} \mathrm{d}t^{'} + e^{-t/\tau} \boldsymbol{f}_{k0} (\boldsymbol{x}-\boldsymbol{u}_k t).
\label{eq:integral_solution}
\end{equation}
Here, $\boldsymbol{f}_{0k}$ is the initial non-equilibrium distribution at $t=0$, and $\boldsymbol{g}_{0k}(\boldsymbol{x}^{'},t^{'})$ is the equilibrium distribution along the characteristic line $\boldsymbol{x}^{'} = \boldsymbol{x}-\boldsymbol{u}_kt^{'}$. Expand the equilibrium and non-equilibrium distribution function by the Taylor series to the second-order accuracy, we have
$$
\boldsymbol{g}_{0k}(\boldsymbol{x}^{'},t^{'}) = \boldsymbol{g}_{0k}(\boldsymbol{x},t) + \nabla\boldsymbol{g}_{0k} \cdot (\boldsymbol{x}^{'}-\boldsymbol{x}) + \frac{\partial \boldsymbol{g}_{0k}}{\partial t}(t^{'} - t),
$$
and 
$$
\boldsymbol{f}_{0k}(\boldsymbol{x}-\boldsymbol{u}_k t) = \boldsymbol{f}_{0k}(\boldsymbol{x}) - \nabla\boldsymbol{f} \cdot \boldsymbol{u}_kt
= \left\{
\begin{array}{ll}
  \boldsymbol{g}_k^l-\nabla \boldsymbol{g}_k^l \cdot \boldsymbol{u}_kt,   &  \boldsymbol{u}_k\cdot \boldsymbol{n}_s \leq 0\\
  \boldsymbol{g}_k^r-\nabla \boldsymbol{g}_k^r \cdot \boldsymbol{u}_kt,   &  \boldsymbol{u}_k\cdot \boldsymbol{n}_s < 0\\
\end{array}.
\right.
$$
The equilibrium distribution $\boldsymbol{g}_{0k}$ at the cell interface is obtained by conservation constraint

$$
<\boldsymbol{g}_{0k}> = <\boldsymbol{g}_k^l>_{>0} + <\boldsymbol{g}_k^r>_{<0},
$$
where $\boldsymbol{g}^l$ and $\boldsymbol{g}^r$ are the equilibrium distributions on the left and right side of the cell interface. The equilibrium spatial derivative on the cell interface is calculated as
$$
<\nabla\boldsymbol{g}_{0k}> =<\nabla\boldsymbol{g}_k^l>_{>0} + <\nabla\boldsymbol{g}_k^r>_{<0}.
$$
The spatial derivative of $\boldsymbol{g}^{l}$ and $\boldsymbol{g}^r$ can be obtained
$$
\nabla\boldsymbol{g}_{k}^{l,r} = \nabla\left[\frac{\boldsymbol{U}_k^{l,r}}{4} + \frac{1}{4c^2}\mathbbm{F}(\boldsymbol{U}^{l,r})\cdot\boldsymbol{u}_k\right],\quad
$$
The right-hand side of the above function requires the calculation of $\nabla\boldsymbol{U}^{l,r}$, and they are calculated by the Green-Gauss method or other high-order gradient schemes. The compatibility condition can obtain the temporal derivatives
$$
<\frac{\partial \boldsymbol{g}_{0k}}{\partial t} + \boldsymbol{u}_k \cdot \nabla \boldsymbol{g}_{0k}> = 0.
$$
Substituting the reconstructed equilibrium and non-equilibrium distribution into Equation \ref{eq:integral_solution}, we have 
\begin{equation}
\begin{aligned}
 \boldsymbol{f}_k(0,t) &= c_1 \boldsymbol{g}_{0k} + c_2 \boldsymbol{u}_k\cdot\nabla\boldsymbol{g}_{0k} + c_3 \frac{\boldsymbol{g}_{0k}}{\partial t} \\&\quad+ c_4 \left[H(u)\boldsymbol{g}^l + (1-H(u)) \boldsymbol{g}^r\right] - c_5 \left[H(u)(\nabla \boldsymbol{g}^l \cdot \boldsymbol{u}_k t + (1-H(u)) \nabla \boldsymbol{g}^r \cdot \boldsymbol{u}_k t\right].
\end{aligned}
\label{eq:integral_solution_taylor}
\end{equation}
The coefficients are given as 
$$
\begin{aligned}
    &c_1 = 1-e^{-t/\tau} ,\\
    &c_2 = (\tau +t)e^{-t/\tau} - \tau ,\\
    &c_3 = t + \tau(e^{-t/\tau}-1),\\
    &c_4 =  e^{t/\tau} ,\\
    &c_5 = te^{t/\tau} .\\
\end{aligned}
$$
Then the numerical flux on the interface can be evaluated as
\begin{equation}
\begin{aligned}
\mathscr{F}_{s}
&= \frac{1}{\Delta t}\int_{t^n}^{t^{n+1}}<\boldsymbol{u}_k\cdot \boldsymbol{n}_s \boldsymbol{f}_k(0,t)>\mathrm{dt}\\
&= \frac{1}{\Delta t}\int_0^{\Delta t}<\boldsymbol{u}_k\cdot \boldsymbol{n}_s \boldsymbol{f}_k(0,t)>\mathrm{dt}\\
&= \frac{1}{\Delta t}[C_1 <u_k\boldsymbol{g}_{0k}> + C_2 <u_k\boldsymbol{u}_k\cdot\nabla \boldsymbol{g}_{0k}>+C_3<u_k\boldsymbol{g}_t>\\ &\quad+ C_4 (<u_k\boldsymbol{g}^l>_{>0}+<u_k\boldsymbol{g}^r>_{<0}) - C_5(<u_k\boldsymbol{u}_k\cdot\nabla\boldsymbol{g}^l>_{>0}+<u_k\boldsymbol{u}_k\cdot\nabla\boldsymbol{g}^r>_{<0})].\\
\end{aligned}
\label{eq:numerical_flux_solution}
\end{equation}
The coefficients are
$$
\begin{aligned}
    &C_1 = \int_0^{\Delta t} (1-e^{-t/\tau}) dt = \Delta t - \tau(1-e^{-\Delta t/\tau}),\\
    &C_2 = \int_0^{\Delta t} [(\tau +t)e^{-t/\tau} - \tau ]dt = \tau(-\Delta t + 2\tau (1-e^{-\Delta t/\tau})-\Delta t e^{-\Delta t/\tau}),\\
    &C_3 = \int_0^{\Delta t} [t + \tau(e^{-t/\tau}-1)] dt = \frac{1}{2}\Delta t^2 - \tau \Delta t - \tau^2 e^{-\Delta t / \tau} + \tau^2,\\
    &C_4 = \int_0^{\Delta t} e^{t/\tau} dt =   \tau(1-e^{-\Delta t/\tau}),\\
    &C_5 = \int_0^{\Delta t} te^{t/\tau} dt = -\tau \Delta t e^{-\Delta t / \tau} + \tau^2(1-e^{-\Delta t / \tau}).
\end{aligned}
$$
We have now solved the numerical flux at the cell interface. 

For the implementation of a non-reflecting boundary condition, the distribution function corresponding to incoming kinetic velocities is set to zero within the boundary ghost cells. This ensures that no information propagates from the boundary into the computational domain. This approach, inherent to the kinetic formulation, simplifies implementation considerably when compared to methods such as FDTD, which typically require additional Perfectly Matched Layer (PML) absorption regions to achieve a similar non-reflecting effect.

\subsection{Lattice Boltzmann method}

The operator-splitting method is employed in LBM to solve the kinetic model Eq.\eqref{eq:discretized vector BGK},
\begin{equation}
\left\{
\begin{aligned}
    \mathcal{R}(\Delta t): &\quad \frac{\partial \boldsymbol{f}_k}{\partial t} = \frac{\boldsymbol{g}_k -\boldsymbol{f}_k}{\tau}, \\
    \mathcal{T}(\Delta t): &\quad \frac{\partial \boldsymbol{f}_k}{\partial t} + \boldsymbol{u}_k\cdot\nabla\boldsymbol{f}_k = 0,
\end{aligned}
\right.
\end{equation}
where $\mathcal{R}$ means collision relaxation operator, which is a local homogeneous relaxation term. $\mathcal{T}$ means transport operator. For the collision operator $\mathcal{R}$, the Crank-Nicolson scheme is employed.
\begin{equation}
    \boldsymbol{f}_{k,i}^{*} =  \frac{1-\frac{\Delta t}{2\tau}}{1+\frac{\Delta t}{2\tau}} \boldsymbol{f}_{k,i}^{n} + \frac{\frac{2\Delta t}{\tau}}{1+\frac{\Delta t}{2\tau}}\boldsymbol{g}_{k,i}^{n}\xrightarrow{\omega=\frac{\frac{2\Delta t}{\tau}}{1+\frac{\Delta t}{2\tau}}} (1-\omega)\boldsymbol{f}_{k,i}^{n}+\omega \boldsymbol{g}_{k,i}^{n},
\end{equation}
where $i$ is the cell index. For the transport operator $\mathcal{T}$, the exact solution can be utilized on structured meshes. Since the transport equation is linear, the analytical solution based on the characteristic theory can be obtained. The solution is
\begin{equation}
    \boldsymbol{f}_k(\boldsymbol{x},t+\Delta t) = \boldsymbol{f}_k(\boldsymbol{x}-\boldsymbol{u}_k \Delta t,t).
    \label{eq:exact transport}
\end{equation}
Let $\Delta \boldsymbol{x} / \Delta t = |\boldsymbol{u}_k|$, Eq.\eqref{eq:exact transport} can be written as
\begin{equation}
    \boldsymbol{f}_k(\boldsymbol{x},t+\Delta t) = \boldsymbol{f}_k(\boldsymbol{x}-\boldsymbol{\hat{e}}_k \Delta x,t),
    \label{eq:lbm transport}
\end{equation}
where $\boldsymbol{\hat{e}}_k$ is the unit direction vector of microscopic velocity $\boldsymbol{u}_k$.
Section \ref{sec:asymptotic analysis} shows that when $\omega=2$, the LBM is second-order in time and space. Additionally, the same analysis indicates that LBM requires the magnitude of its microscopic velocities to exceed the speed of light $c$ to ensure stability.
Solving the transport equation exactly realizes the multidimensionality of the scheme. However, this method is not readily generalizable to unstructured meshes.

After the transport step, map the microscopic distribution function to the macroscopic quantities by the following formula.
$$
\boldsymbol{U}_i^{*} = \sum_k \boldsymbol{f}_k.
$$

\subsection{Flux vector splitting scheme}
\label{sec:flux vector splitting scheme}

In contrast to GKS, the flux vector splitting (FVS) scheme uses a direct eigendecomposition to split the flux vector. The governing one-dimensional system is given by
\begin{equation}
\frac{\partial \boldsymbol{U}}{\partial t} + \mathbbm{A}_1\frac{\partial \boldsymbol{U}}{\partial x} = 0,
\label{eq:PHM1D}
\end{equation}
where
\begin{equation*}
    \mathbbm{A}_{1}=\left(\begin{array}{cccccccc}
    0 & 0 & 0 & 0 & 0 & 0 & c^{2} \chi & 0 \\
    0 & 0 & 0 & 0 & 0 & c^{2} & 0 & 0 \\
    0 & 0 & 0 & 0 & -c^{2} & 0 & 0 & 0 \\
    0 & 0 & 0 & 0 & 0 & 0 & 0 & \gamma \\
    0 & 0 & -1 & 0 & 0 & 0 & 0 & 0 \\
    0 & 1 & 0 & 0 & 0 & 0 & 0 & 0 \\
    \chi & 0 & 0 & 0 & 0 & 0 & 0 & 0 \\
    0 & 0 & 0 & c^{2} \gamma & 0 & 0 & 0 & 0
    \end{array}\right).
\end{equation*}
Since Eq. \eqref{eq:PHM1D} is linear, its Riemann solution can be captured exactly without iteration. The numerical flux at an interface is computed using an upwind scheme,
\begin{equation*}
    \mathscr{F}_s(\boldsymbol{U}_L, \boldsymbol{U}_R) = \frac{1}{2}\mathbbm{A}_1(\boldsymbol{U}_L + \boldsymbol{U}_R) - \frac{1}{2}|\mathbbm{A}_1|(\boldsymbol{U}_R - \boldsymbol{U}_L),
\end{equation*}
where $|\mathbbm{A}_1| = \mathbbm{A}_1^{+} - \mathbbm{A}_1^{-}$. The matrices $\mathbbm{A}_1^{\pm}$ are constructed via $\mathbbm{A}_1^{\pm} = \mathbbm{R}\Lambda^{\pm}\mathbbm{R}^{-1}$, with $\Lambda^{\pm} = \text{diag}(\lambda_1^{\pm},\cdots,\lambda_8^{\pm})$. Here, $\lambda_m$ is the $m$-th eigenvalue of $\mathbbm{A}_1$, $\lambda_m^{\pm}=(\lambda_m\pm|\lambda_m|)/2$, and $\mathbbm{R}=[\boldsymbol{r}_1,\cdots, \boldsymbol{r}_8]$ is the right eigenvector matrix. Second-order spatial accuracy is attained by computing the average gradient within each control volume using either the Green-Gauss or the Least-Squares method. These cell gradients are subsequently employed in a piecewise linear reconstruction to determine the solution states at cell interfaces. For temporal integration, the second-order accurate Runge-Kutta method is applied as the time-stepping scheme.

For non-reflecting boundary conditions, incoming characteristic information is suppressed. At a domain boundary interface, let \(\boldsymbol{U}_R\) denote the ghost cell state outside the computational domain and \(\boldsymbol{U}_L\) the interior cell state adjacent to the boundary. The numerical flux at this interface is then computed as
\[
\mathscr{F}_s = \mathbb{A}^{+} \boldsymbol{U}_L,
\]
where \(\mathbb{A}^{+}\) contains only the non-negative eigenvalues of the flux Jacobian. This formulation ensures that outgoing-wave information is transmitted while incoming-wave information is eliminated.

\section{Analysis}
\label{sec:analysis}
\subsection{Asymptotic analysis of GKS}
\label{sec:analysis gks}

This section provides an analysis of the proposed GKS formulation under different regimes of the relaxation time \( \tau \).

\textbf{Non-viscous limit (\( \tau = 0 \))}: In the limit of zero relaxation time, the exponential term \( e^{-\Delta t/\tau} \) vanishes. Consequently, the free-transport component of the flux in Eq. \eqref{eq:integral_solution} disappears, leaving only the equilibrium contribution. As shown by Eq. \eqref{eq:numerical_flux_solution}, the numerical flux simplifies to
$$
\mathscr{F}_{s} = \frac{1}{\Delta t} \left\langle u_k \boldsymbol{g}_{0k} + \frac{\Delta t}{2} \frac{\partial \boldsymbol{g}_{0k}}{\partial t} \right\rangle,
$$
which corresponds to the numerical flux in the non-viscous, continuum limit.

\textbf{Viscous limit (small, finite \( \tau \))}: For small but finite relaxation times, the free-transport flux decays exponentially. The equilibrium component thus dominates the numerical flux. In Eq. \eqref{eq:numerical_flux_solution}, the coefficients \( C_1 \), \( C_2 \), and \( C_3 \) are dominated by \( \Delta t \), \( -\tau \Delta t \), and \( \frac{\Delta t^2}{2} - \tau \Delta t \), respectively. This yields the equilibrium flux
$$
\mathscr{F}_{s} = \frac{1}{\Delta t}  \left\langle \boldsymbol{g}_{0k} -\tau \left(\boldsymbol{u}\cdot\nabla\boldsymbol{g}_{0k} + \frac{\partial \boldsymbol{g}_{0k}}{\partial t} \right) + \frac{\Delta t}{2} \frac{\partial \boldsymbol{g}_{0k}}{\partial t} \right\rangle,
$$
which is consistent with the numerical flux derived for the viscous regime \cite{xu2001}.

\textbf{Collisionless regime (large \( \tau \))}: In the limit of very large relaxation times, \( e^{-\Delta t/\tau} \to 1 \). The free-transport term in Eq. \eqref{eq:integral_solution} therefore dominates, and the GKS reduces to a kinetic vector flux splitting (KFVS) scheme. The following discussion establishes the equivalence between a KFVS scheme in a D1Q2 discrete-velocity space and the standard flux-vector splitting (FVS) method for the one-dimensional Maxwell equations. Consider the one-dimensional transverse electromagnetic mode described by:
\[
\begin{aligned}
&\frac{\partial E_y}{\partial t} + c^2 \frac{\partial B_z}{\partial x} = 0, \\
&\frac{\partial B_z}{\partial t} + \frac{\partial E_y}{\partial x} = 0.
\end{aligned}
\]
The velocity space is discretized using two points:
\[
\boldsymbol{u}_1 = (c,0,0)^T, \quad \boldsymbol{u}_2 = (-c,0,0)^T.
\]
The macroscopic variables are mapped to microscopic distribution functions via:
\[
\begin{aligned}
&g_{E_y,1} = \frac{E_y}{2} + \frac{c B_z}{2}, g_{E_y,2} = \frac{E_y}{2} - \frac{c B_z}{2}, \\
&g_{B_z,1} = \frac{B_z}{2} + \frac{E_y}{2c}, g_{B_z,2} = \frac{B_z}{2} - \frac{E_y}{2c},
\end{aligned}
\]
where subscripts 1 and 2 denote distributions corresponding to \(\boldsymbol{u}_1\) and \(\boldsymbol{u}_2\) respectively.
This decomposition yields left-going fluxes across an interface:
\[
\sum_{\boldsymbol{u}_k>0} \boldsymbol{u}_k g_{E_y,k} = \frac{cE_y}{2} + \frac{c^2 B_z}{2}, \quad
\sum_{\boldsymbol{u}_k>0} \boldsymbol{u}_k g_{B_z,k} = \frac{cB_z}{2} + \frac{E_y}{2},
\]
and right-going fluxes:
\[
\sum_{\boldsymbol{u}_k<0} \boldsymbol{u}_k g_{E_y,k} = -\frac{cE_y}{2} + \frac{c^2 B_z}{2}, \quad
\sum_{\boldsymbol{u}_k<0} \boldsymbol{u}_k g_{B_z,k} = -\frac{cB_z}{2} + \frac{E_y}{2}.
\]
On the other hand, the system Jacobian is:
\[
\mathbb{A} = \begin{pmatrix}
0 & c^2 \\
1 & 0
\end{pmatrix},
\]
with eigendecomposition:
\[
\Lambda = \begin{pmatrix}
c & 0 \\
0 & -c
\end{pmatrix}, \quad
\mathbb{R} = \begin{pmatrix}
c & -c \\
1 & 1
\end{pmatrix}, \quad
\mathbb{R}^{-1} = \frac{1}{2c} \begin{pmatrix}
1 & c \\
-1 & c
\end{pmatrix}.
\]
The split Jacobians become:
\[
\mathbb{A}^+ = \mathbb{R} \Lambda^+ \mathbb{R}^{-1} = \frac{1}{2} \begin{pmatrix}
c & c^2 \\
1 & c
\end{pmatrix}, \quad
\mathbb{A}^- = \mathbb{R} \Lambda^- \mathbb{R}^{-1} = \frac{1}{2} \begin{pmatrix}
-c & c^2 \\
1 & -c
\end{pmatrix}.
\]
The corresponding wave decompositions are:
\[
\mathbb{A}^+ \boldsymbol{U} = \frac{1}{2} \begin{pmatrix}
c & c^2 \\
1 & c
\end{pmatrix} \begin{pmatrix}
E_y \\
B_z
\end{pmatrix} = \begin{pmatrix}
\displaystyle \frac{cE_y}{2} + \frac{c^2 B_z}{2} \\[6pt]
\displaystyle \frac{E_y}{2} + \frac{c B_z}{2}
\end{pmatrix},
\]
and
\[
\mathbb{A}^- \boldsymbol{U} = \frac{1}{2} \begin{pmatrix}
-c & c^2 \\
1 & -c
\end{pmatrix} \begin{pmatrix}
E_y \\
B_z
\end{pmatrix} = \begin{pmatrix}
\displaystyle -\frac{cE_y}{2} + \frac{c^2 B_z}{2} \\[6pt]
\displaystyle \frac{E_y}{2} - \frac{c B_z}{2}
\end{pmatrix}.
\]
These expressions exactly match the kinetic fluxes derived above, demonstrating the equivalence between KFVS and FVS for this one-dimensional system. In multidimensional configurations, the kinetic approach provides a more general framework for flux decomposition.

\subsection{Asymptotic analysis of LBM}
\label{sec:asymptotic analysis}

In this section, we also employ the one-dimensional polarized pure electromagnetic wave equation in section \ref{sec:analysis gks} to illustrate the asymptotic behavior of the LBM method.
Consider the kinetic equation of $E_y$ first,
\begin{align}
    &\frac{\partial f_{E_y,1}}{\partial t} -\lambda c \frac{\partial f_{E_y,1}}{\partial x} = \frac{g_{E_y,1} - f_{E_y,1}}{\tau}, \label{eq:EyV1}\\
    &\frac{\partial f_{E_y,2}}{\partial t} +\lambda c \frac{\partial f_{E_y,2}}{\partial x} = \frac{g_{E_y,2} - f_{E_y,2}}{\tau} \label{eq:EyV2},
\end{align}
where $\lambda$ is a parameter used to adjust the microscopic velocity for stability of LBM. 
According to Eq.\eqref{eq:def_f_discretized}, $E_y=f_{E_y,1}+f_{E_y,2}$, define $q = \lambda c f_{E_y,1} - \lambda c f_{E_y,2}$, then Eq.\eqref{eq:EyV1} + Eq.\eqref{eq:EyV2} and (Eq.\eqref{eq:EyV1} - Eq.\eqref{eq:EyV2} )$\times \lambda c$ give the following equations,
\begin{equation}
    \frac{\partial E_y}{\partial t} -\frac{\partial q}{\partial x} = 0,\quad
    \frac{\partial q}{\partial t} -\lambda^2 c^2 \frac{\partial E_y}{\partial x} = \frac{c^2B_z - q}{\tau}.
    \label{eq:1Dequiv}
\end{equation}
Expand
$$
q = q_0 + \tau q_1 + \mathcal{O}(\tau^2).
$$
Substitute the above expansion into Eq.\eqref{eq:1Dequiv}, for $\mathcal{O}(1)$ balance, we have
$$
q_0 = c^2 B_z.
$$
For $\mathcal{O}(\tau)$ balance, we have
$$
q_1 = - \frac{\partial q_0}{\partial t} + \lambda^2 c^2 \frac{\partial E_y}{\partial x} = \lambda^2 c^2 \frac{\partial E_y}{\partial x} - \frac{\partial c^2 B_z}{\partial t} = (\lambda^2-1)c^2\frac{\partial E_y}{\partial x}.
$$
Therefore, the first order approximation is $q=c^2B_z -\tau (\lambda^2-1)c^2\frac{\partial E_y}{\partial x} + \mathcal{O}(\tau^2)$, and the macroscopic equation of $E_y$ becomes
$$
\frac{\partial E_y}{\partial t} +\frac{\partial c^2 B_z}{\partial x} =\tau\frac{\partial}{\partial x}\left(\left(\lambda^2- 1\right)c^2 \frac{\partial E_y}{\partial x} \right).
$$
Note that $\lambda > c$ is required for the kinetic approximation to be dissipative, which is referred to as the sub-characteristic condition by Liu \cite{liu1987hyperbolic}. This condition ensures that the relaxation system remains stable and dissipative. If $\tau = 0$ or $\lambda = 1$, the dissipation term vanishes, and the system theoretically recovers the original Maxwell equations accurately. 

Next, we analyze the numerical error following the analysis in \cite{graille_approximation_2014}. The first step is a relaxation step,
\begin{equation}
f_{E_y,k}^{*}(x_i, t) = (1-\omega)f_{E_y, k}(x_i, t) + \omega g_{E_y, k}(x_i , t) ,
\label{eq:err_relax}
\end{equation}
then the Maxwell step is transport, for the sake of simplicity, we suppose this step is exact, i.e.
\begin{equation}
f_{E_y,k}(x_i, t+\Delta t) = f_{E_y, k}^{*}(x_i - u_k \Delta t, t).
\label{eq:err_transport}
\end{equation}
Expand Eq.\eqref{eq:err_transport} with respect to $\Delta t$,
\begin{equation}
f_{E_y,k} + \frac{\partial f_{E_y,k}}{\partial t}\Delta t + \frac{1}{2}\frac{\partial^2 f_{E_y,k}}{\partial t^2}\Delta t^2 = f_{E_y,k}^{*} - \frac{\partial f_{E_y,k}^{*}}{\partial x}u_k\Delta t + \frac{1}{2}\frac{\partial^2 f_{E_y,k}^{*}}{\partial x^2}u_k^2\Delta t^2 + \mathcal{O}(\Delta t^3).
\label{eq:err_taylor}
\end{equation}
Take the zeroth order moment, we have
\begin{equation}
E_y + \frac{\partial E_y}{\partial t}\Delta t + \frac{1}{2}\frac{\partial^2 E_y}{\partial t^2}\Delta t^2 = E_y^{*} - \frac{\partial q^{*}}{\partial x}\Delta t + \frac{1}{2}\frac{\partial^2 E_y^{*}}{\partial x^2}\lambda^2c^2\Delta t^2 + \mathcal{O}(\Delta t^3).
\label{eq:err_taylor_0th_moment}
\end{equation}
Take the first order moment, we have
\begin{equation}
q + \frac{\partial q}{\partial t}\Delta t + \frac{1}{2}\frac{\partial^2 q}{\partial t^2}\Delta t^2 = q^{*} - \frac{\partial E_y^{*}}{\partial x}\lambda^2 c^2\Delta t + \frac{1}{2}\frac{\partial^2 q^{*}}{\partial x^2}\lambda^2c^2\Delta t^2 + \mathcal{O}(\Delta t^3).
\label{eq:err_taylor_1th_moment}
\end{equation}
For $\mathcal{O}(\Delta t)$ order of accuracy, Eq.\eqref{eq:err_taylor_0th_moment} gives $E_y = E_y^{*}$ which is obvious since the relaxation process \eqref{eq:err_relax} doesn't change $E_y$. Eq.\eqref{eq:err_taylor_1th_moment} gives $q = q^{*} + \mathcal{O}(\Delta t)$. Note that the first order moment of Eq.\eqref{eq:err_relax} gives
\begin{equation}
q^{*} = (1-\omega)q + \omega c^2 B_z,
\label{eq:err_relax_1st_moment}
\end{equation}
therefore
$$q = c^2B_z + \mathcal{O}(\Delta t), \quad q^{*} = c^2B_z + \mathcal{O}(\Delta t).$$
For $\mathcal{O}(\Delta t^2)$ order of accuracy, combining with the above formula, Eq.\eqref{eq:err_taylor_0th_moment} gives
$$
\frac{\partial E_y}{\partial t} + \frac{\partial c^2B_z}{\partial x} = \mathcal{O}(\Delta t).
$$
Eq.\eqref{eq:err_taylor_1th_moment} gives
$$
q - q^{*} = -(\frac{\partial B_z}{\partial t}+\lambda^2\frac{\partial E_y}{\partial x})c^2\Delta t + \mathcal{O}(\Delta t^2) \equiv -\theta \Delta t + \mathcal{O}(\Delta t^2),
$$
where $\theta \equiv (\frac{\partial B_z}{\partial t}+\lambda^2\frac{\partial E_y}{\partial x})c^2$. Substitute Eq.\eqref{eq:err_relax_1st_moment} into the above equation, we have
\begin{equation}
    q = c^2B_z - \frac{\theta\Delta t}{\omega t} + \mathcal{O}(\Delta t^2),\quad q^{*} = c^2B_z - \frac{1-\omega}{\omega}\theta \Delta t + \mathcal{O}(\Delta t^2).
\end{equation}
For $\mathcal{O}(\Delta t^3)$ order of accuracy, Eq.\eqref{eq:err_taylor_0th_moment} gives
\begin{align}
    \begin{split}
        \frac{\partial E_y}{\partial t} + \frac{\partial c^2B_z}{\partial x} &= (\frac{1}{2}\frac{\partial^2 E_y}{\partial x^2}\lambda^2 c^2 - \frac{1}{2}\frac{\partial^2 E_y}{\partial t^2} + \frac{1-\omega}{\omega}\theta)\Delta t+ \mathcal{O}(\Delta t^2)\\
        & = (\frac{1}{2} + \frac{1-\omega}{\omega})\theta\Delta t + \mathcal{O}(\Delta t^2).
    \end{split}
    \label{eq:err_2nd_equiv}
\end{align}
Eq.\eqref{eq:err_2nd_equiv} gives equivalent equation to $\mathcal{O}(\Delta t^2)$. When $\omega=2$, the $\mathcal{O}(\Delta t)$ term vanishes, resulting in a scheme of second-order temporal accuracy.

\section{Numerical tests}
\label{sec:numerical studies}

In the following cases, the following units are used: light speed $ c = 1$ m/s, vacuum permittivity $\epsilon_0 = 1$ F/m, and vacuum permeability $\mu_0 = 1$ H/m. 

\subsection{Plane wave test}
\label{sec:plane wave test}

\begin{figure}
    \begin{center}
    \begin{subfigure}[b]{0.48\textwidth}
 \includegraphics[width=0.9\textwidth]{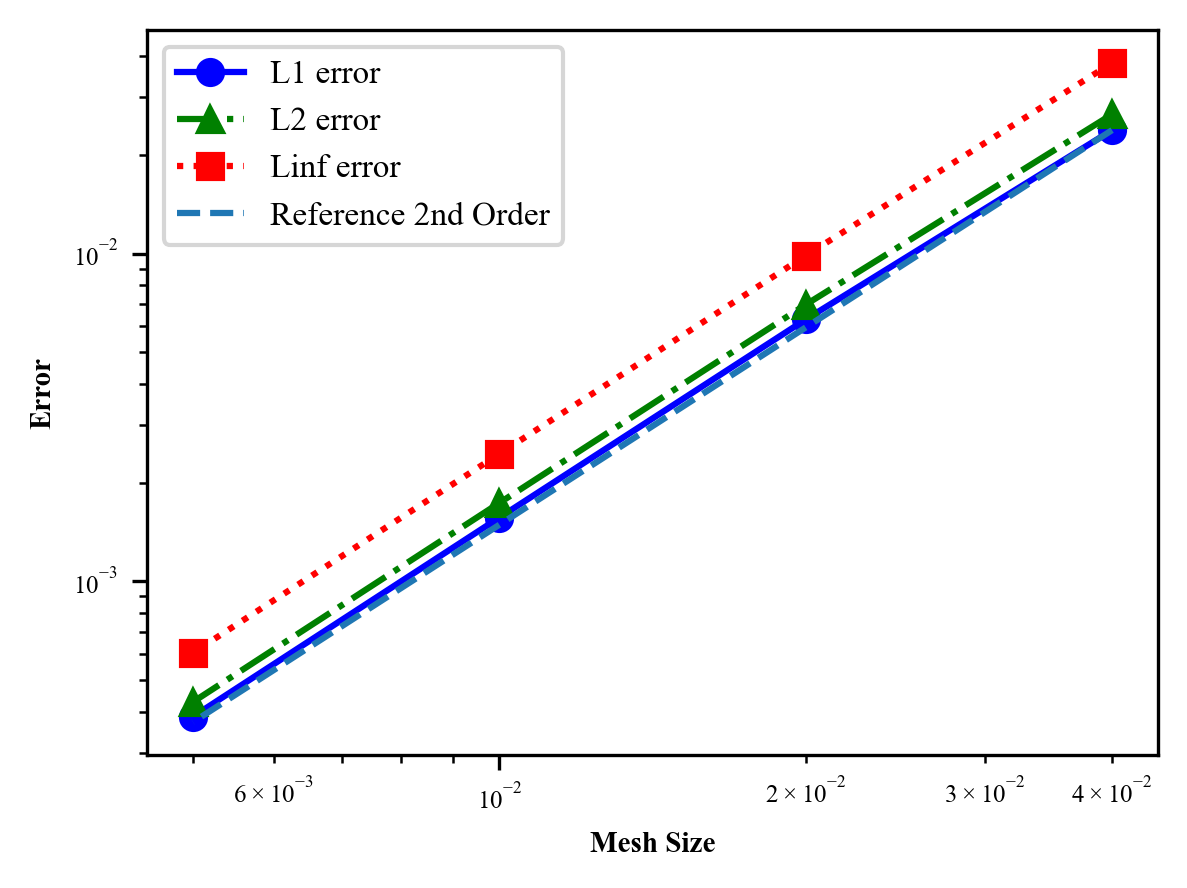}
    \end{subfigure}
    \hfill
    \begin{subfigure}[b]{0.48\textwidth}
 \includegraphics[width=0.9\textwidth]{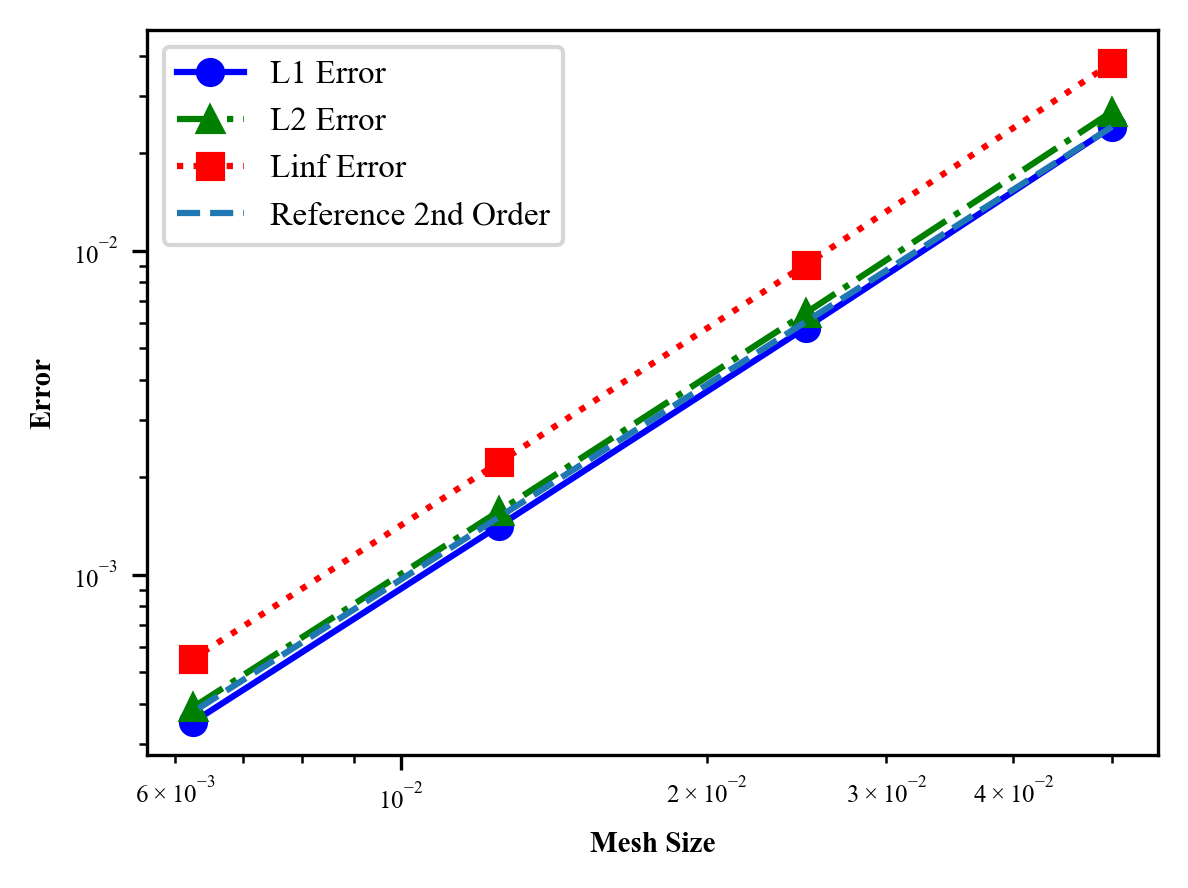}
    \end{subfigure}
    \end{center}
    \caption{Plane wave test: Convergence order of GKS (left) and LBM (right).}
    \label{fig:plane wave order}
\end{figure}

In this section, we use the plane-wave solution of Maxwell's equations on structured meshes to evaluate the convergence order of the GKS and LBM. The spatial domain is defined as $x \in [0,1]$m and $t_{max} = 1$s. We use the following exact time-domain solution for a plane wave:
\begin{equation}
\left\{
\begin{aligned}
    &E_z(\boldsymbol{x},t) = \cos(2\pi(x-ct)),\\
    &B_y(\boldsymbol{x},t) = -\cos(2\pi(x-ct)).
\end{aligned}
\right.
\end{equation}
The CFL number in this test is set as 0.5. 
Mesh number $N=20^3, 40^3, 80^3, 160^3$ is used. The boundary condition is periodic to eliminate boundary effects. Artificial potential parameters are set as $\gamma = 0, \chi = 0$. From Figure \ref{fig:plane wave order} and Table \ref{tab:GKS planewave error} and \ref{tab:LBM planewave error}, we can see that the convergence order of GKS and LBM is second order.

\begin{table}[!h]
    \centering
    \begin{tabular}{cllllll}
    \hline  { Mesh } & $L^1$  { error } &  { Order } & $L^2$  { error } &  { Order } & $L^{\infty}$  { error } & {Order}\\
    \hline 20 & 2.3809$\mathrm{e}{-02}$  &      & 2.6609$\mathrm{e}{-02}$ &      & 3.8171$\mathrm{e}{-02}$ & \\
    40        & 6.2955$\mathrm{e}{-03}$  & 1.92 & 6.9871$\mathrm{e}{-03}$ & 1.93 & 9.8584$\mathrm{e}{-03}$ & 1.95\\
    80        & 1.5651$\mathrm{e}{-03}$  & 2.01 & 1.7346$\mathrm{e}{-03}$ & 2.01 & 2.4483$\mathrm{e}{-03}$ & 2.01\\
    160        & 3.8659$\mathrm{e}{-04}$ & 2.02 & 4.2852$\mathrm{e}{-04}$ & 2.02 & 6.0527$\mathrm{e}{-04}$ & 2.02\\
    \hline
    \end{tabular}
    \caption{Plane wave test: errors and convergence orders at $t=1$ obtained by GKS.}
    \label{tab:GKS planewave error}
\end{table}

\begin{table}[!h]
    \centering
    \begin{tabular}{cllllll}
    \hline  { Mesh } & $L^1$  { error } &  { Order } & $L^2$  { error } &  { Order } & $L^{\infty}$  { error } & {Order}\\
    \hline 20 & 2.4165$\mathrm{e}{-02}$  &      & 2.6915$\mathrm{e}{-02}$ &      & 3.8033$\mathrm{e}{-02}$ & \\
    40        & 5.7943$\mathrm{e}{-03}$  & 2.06 & 6.4295$\mathrm{e}{-03}$ & 2.07 & 9.0707$\mathrm{e}{-03}$ & 2.07\\
    80        & 1.5749$\mathrm{e}{-03}$  & 2.03 & 1.5749$\mathrm{e}{-03}$ & 2.03 & 2.2262$\mathrm{e}{-03}$ & 2.03\\
    160        & 3.9006$\mathrm{e}{-04}$ & 2.01 & 3.9006$\mathrm{e}{-04}$ & 2.01 & 5.5163$\mathrm{e}{-04}$ & 2.01\\
    \hline
    \end{tabular}
    \caption{Plane wave test: errors and convergence orders at $t=1$ obtained by LBM.}
    \label{tab:LBM planewave error}
\end{table}

\subsection{Two-dimensional transverse magnetic modes in a cavity}

In this section, we present a test of the two-dimensional (2D) transverse magnetic (TM) eigenmodes within a rectangular metal cavity. The domain is the unit square \( [0, 1]m \times [0, 1]m \). The electric field is initialized as
\[
E_z = E_0 \sin(a x) \sin(a y),
\]
with \( a = 8\pi \) and \( E_0 = 1 \)V/m. All other field components are zero. The simulation is advanced to \( t = 1.775 \, \text{s} \). The analytical solution for \( E_z \) is given by
\[
E_z = E_0 \sin(a x) \sin(a y) \cos(\omega_0 t),
\]
where the angular frequency is \( \omega_0 = c \sqrt{2a^2} = 8\sqrt{2}\pi \) rad/s.

\begin{figure}
\centering
\begin{subfigure}[b]{0.48\textwidth}
    \includegraphics[width=0.9\textwidth]{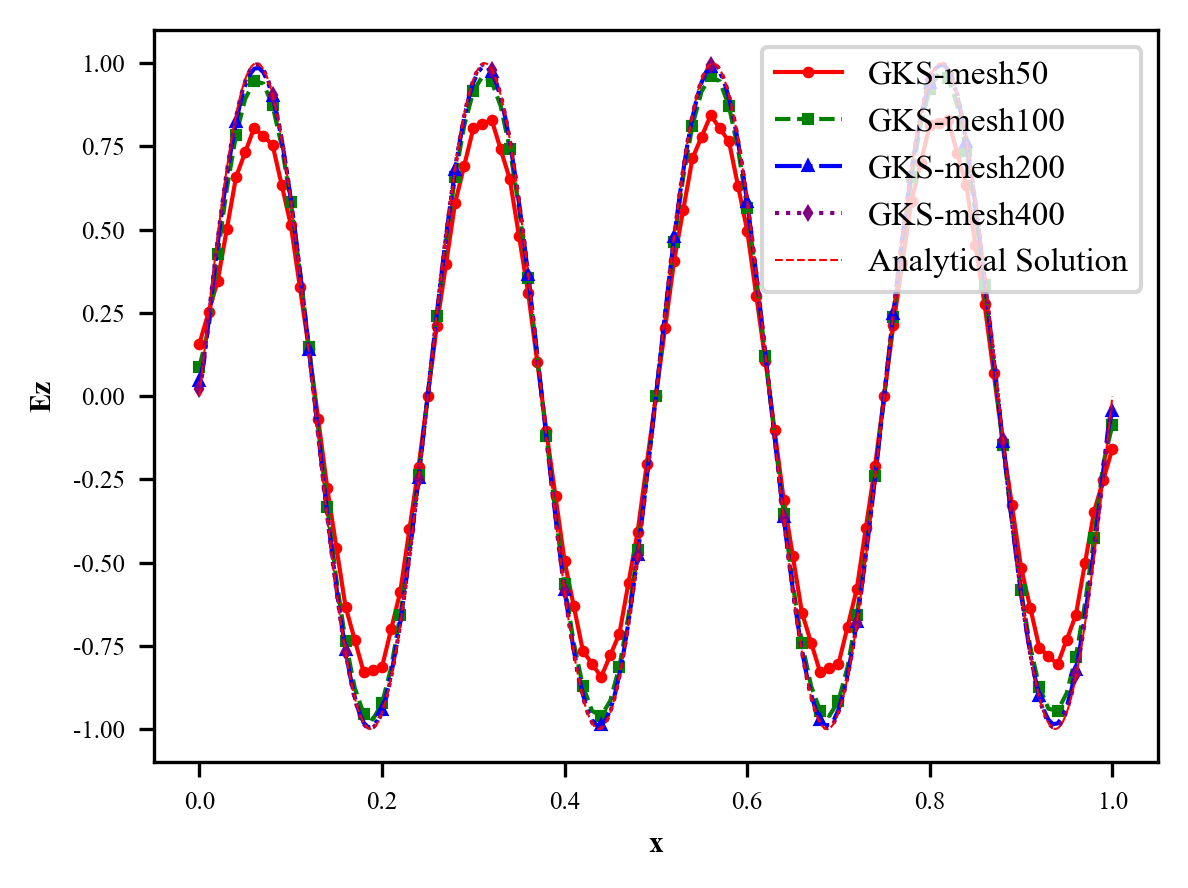} 
\end{subfigure}
\hfill
\begin{subfigure}[b]{0.48\textwidth}
    \includegraphics[width=0.9\textwidth]{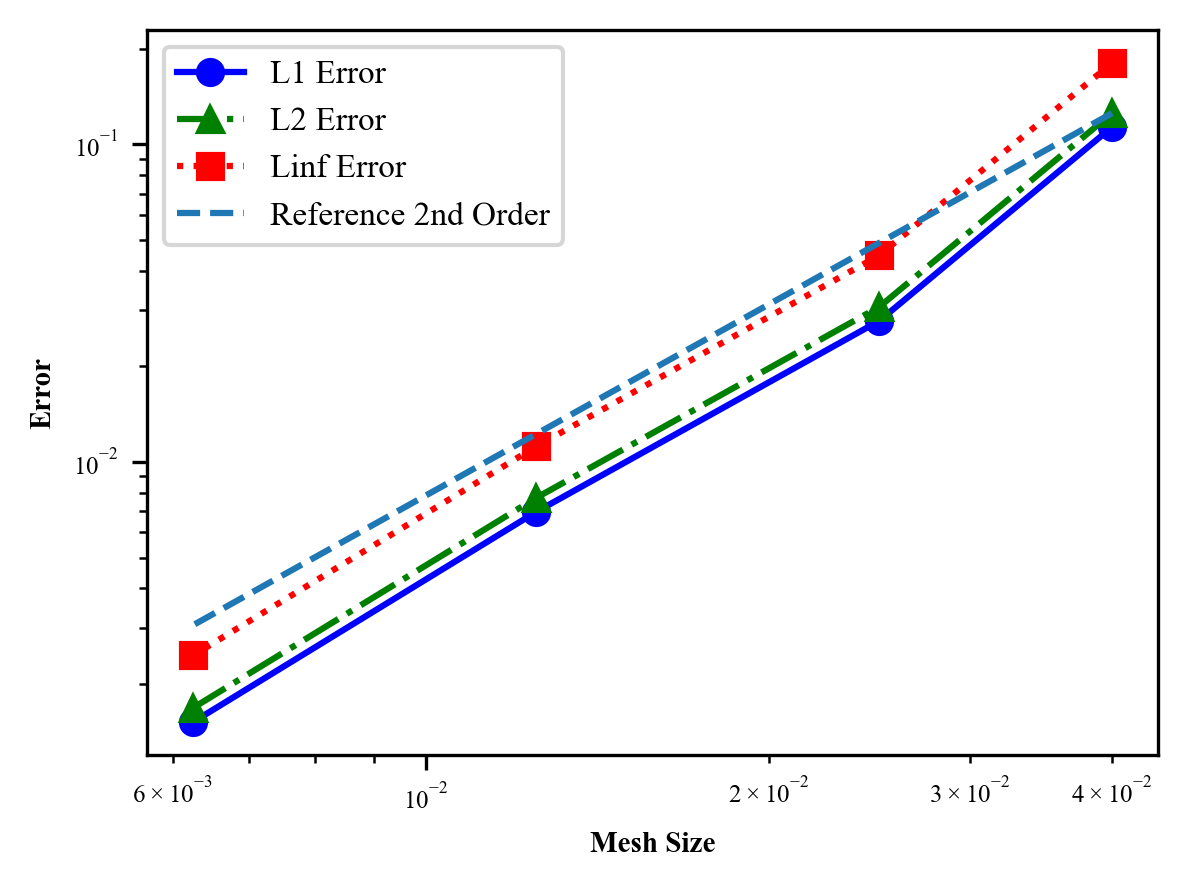} 
\end{subfigure}
\caption{2D TM modes: the $E_z$ profile at $y=0.3125$m (left) and  the convergence order of GKS (right) }
\label{fig:2DTMorder}
\end{figure}

To assess the convergence of the GKS, we perform simulations on successively refined grids: \( 50 \times 50 \), \( 100 \times 100 \), \( 200 \times 200 \), and \( 400 \times 400 \). The time step is selected to satisfy a CFL condition of 0.5, thereby testing the combined accuracy of the spatial and temporal discretization. The results confirm second-order convergence, as summarized in Table \ref{tab:2DTMGKS} and Figure \ref{fig:2DTMorder}.

\begin{figure}
    \centering
    \includegraphics[width=0.7\linewidth]{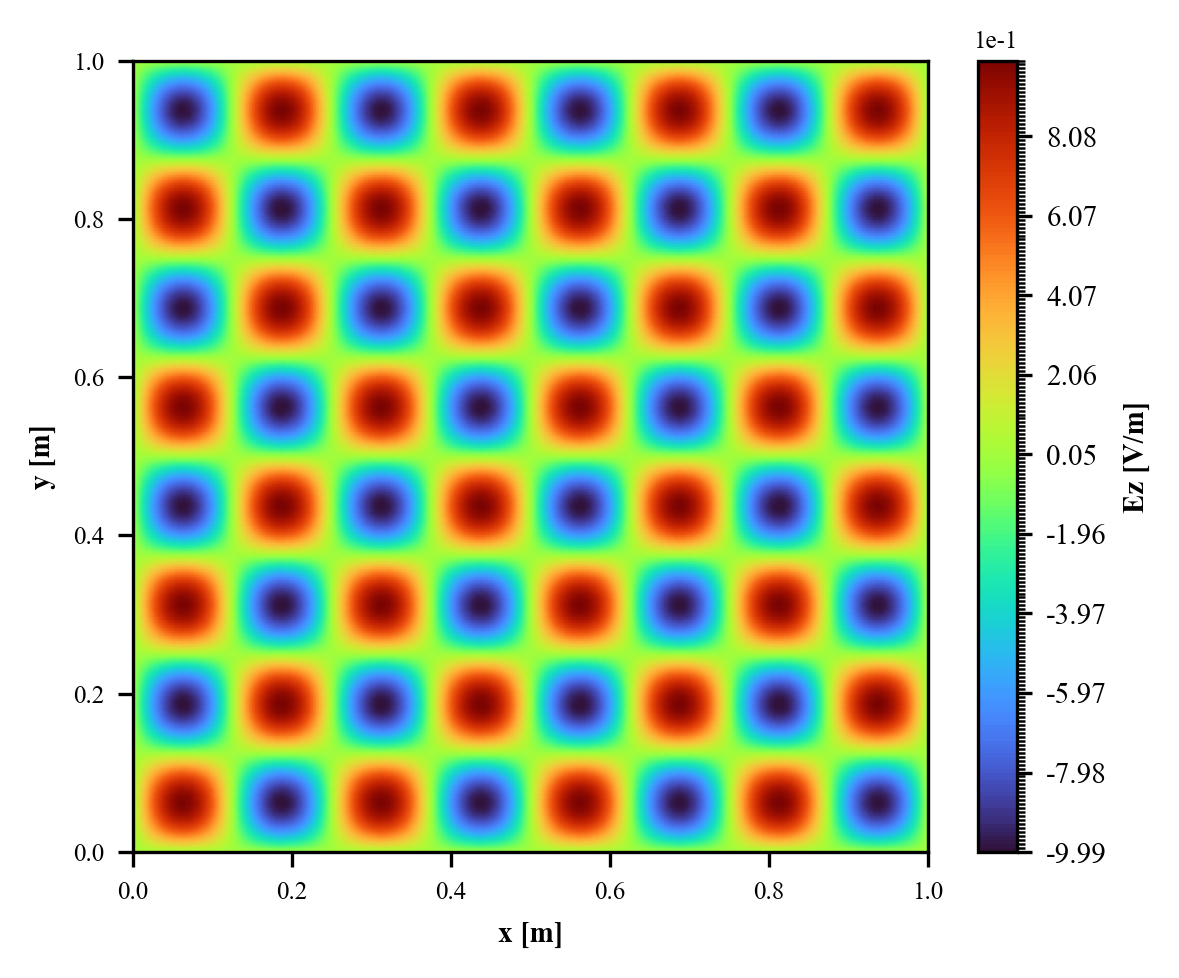}
    \caption{2D TM modes: $E_z$ contour plots at $t=1.775s$.}
    \label{fig:2DTMEzcontour}
\end{figure}

As shown in the left panel of Figure \ref{fig:2DTMorder}, the scheme accurately captures the eigenmode solution even on the coarsest (\( 50 \times 50 \)) mesh. As grid resolution increases, which provides more points per wavelength, the numerical solution converges to the analytical result. The spatial structure of the eigenmode is clearly visible in the contour plot of \( E_z \) at \( t = 1.775 \, \text{s} \), presented in Figure \ref{fig:2DTMEzcontour}.

\begin{table}[!h]
    \centering
    \begin{tabular}{cllllll}
    \hline  { Mesh } & $L^1$  { error } &  { Order } & $L^2$  { error } &  { Order } & $L^{\infty}$  { error } & {Order}\\
    \hline 50 & 1.1341$\mathrm{e}{-01}$  &       & 1.2538$\mathrm{e}{-01}$ &      & 1.8044$\mathrm{e}{-01}$ & \\
    100        & 2.7825$\mathrm{e}{-02}$  & 2.02 & 3.0794$\mathrm{e}{-02}$ & 2.02 & 4.4736$\mathrm{e}{-02}$ & 2.01\\
    200        & 6.9623$\mathrm{e}{-02}$  & 2.00 & 7.7077$\mathrm{e}{-03}$ & 2.00 & 1.1190$\mathrm{e}{-02}$ & 2.00\\
    400        & 1.5168$\mathrm{e}{-03}$  & 2.19 & 1.6805$\mathrm{e}{-03}$ & 2.19 & 2.4592$\mathrm{e}{-03}$ & 2.18\\
    \hline
    \end{tabular}
    \caption{2D TM modes: errors and convergence orders of GKS.}
    \label{tab:2DTMGKS}
\end{table}

\subsection{Riemann problem}

This section evaluates the numerical methods using a one-dimensional Riemann problem. Discontinuities are ubiquitous in both fluid dynamics and plasma physics; therefore, this test case demonstrates the potential applicability of the proposed schemes to such problems. The domain is \([0, 1] \, \text{m}\) with a discontinuity initialized at \(x = 0.5 \, \text{m}\). Open boundary conditions are applied at both ends. The initial discontinuity gives rise to a series of waves that propagate at or below the speed of light. The initial left and right states are defined as:
\[
\left[\begin{array}{c}
E_x \\
E_y \\
E_z \\
B_x \\
B_y \\
B_z
\end{array}\right]_l=\left[\begin{array}{c}
0.0 \\
1.0 \\
0.0 \\
1.0 \\
-0.75 \\
0.0
\end{array}\right], \quad
\left[\begin{array}{c}
E_x \\
E_y \\
E_z \\
B_x \\
B_y \\
B_z
\end{array}\right]_r=\left[\begin{array}{c}
0.0 \\
-1.0 \\
0.0 \\
1.0 \\
0.75 \\
0.0
\end{array}\right].
\]

\begin{figure}
\centering
\begin{subfigure}[b]{0.48\textwidth}
    \includegraphics[width=0.95\textwidth]{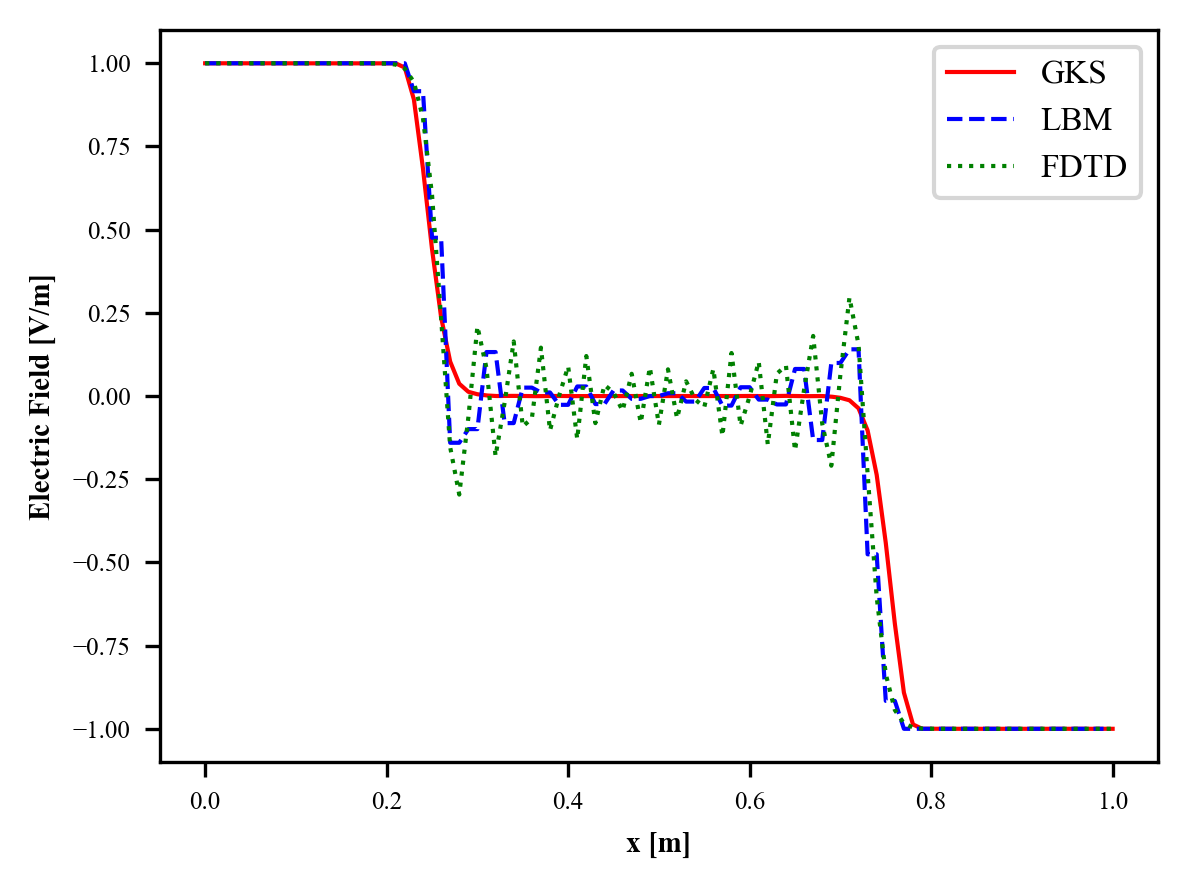} 
    \caption{\(E_y\)}
\end{subfigure}
\hfill
\begin{subfigure}[b]{0.48\textwidth}
    \includegraphics[width=0.95\textwidth]{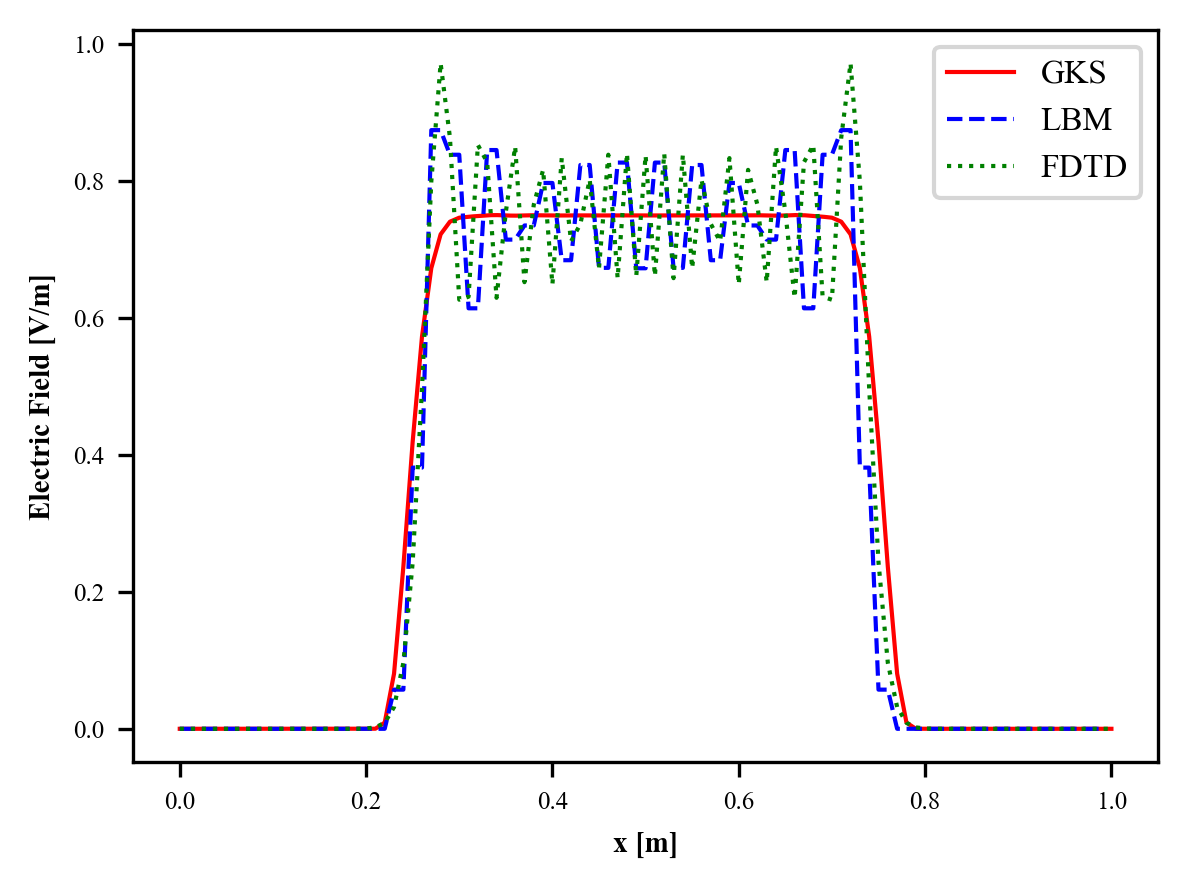} 
    \caption{\(E_z\)}
\end{subfigure}
\vfill
\begin{subfigure}[b]{0.48\textwidth}
    \includegraphics[width=0.95\textwidth]{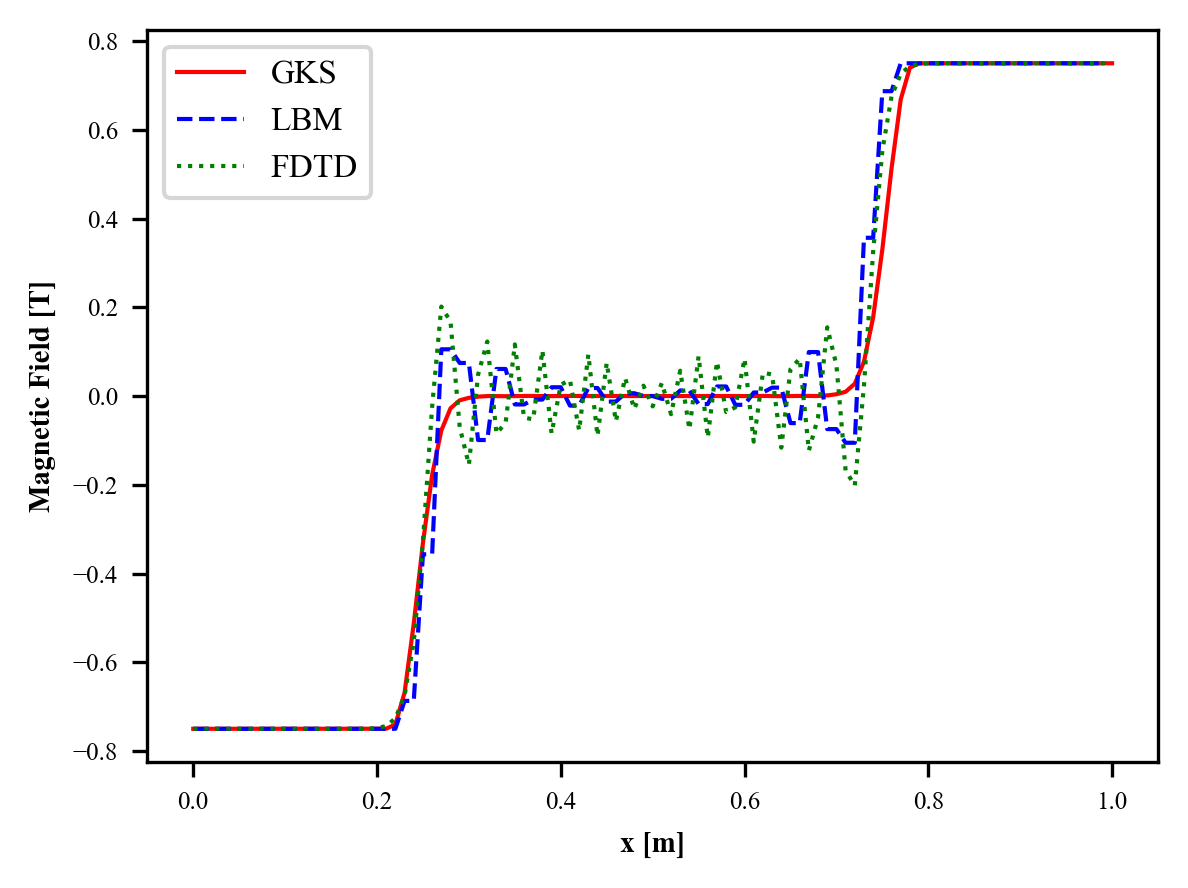} 
    \caption{\(B_y\)}
\end{subfigure}
\hfill
\begin{subfigure}[b]{0.48\textwidth}
    \includegraphics[width=0.95\textwidth]{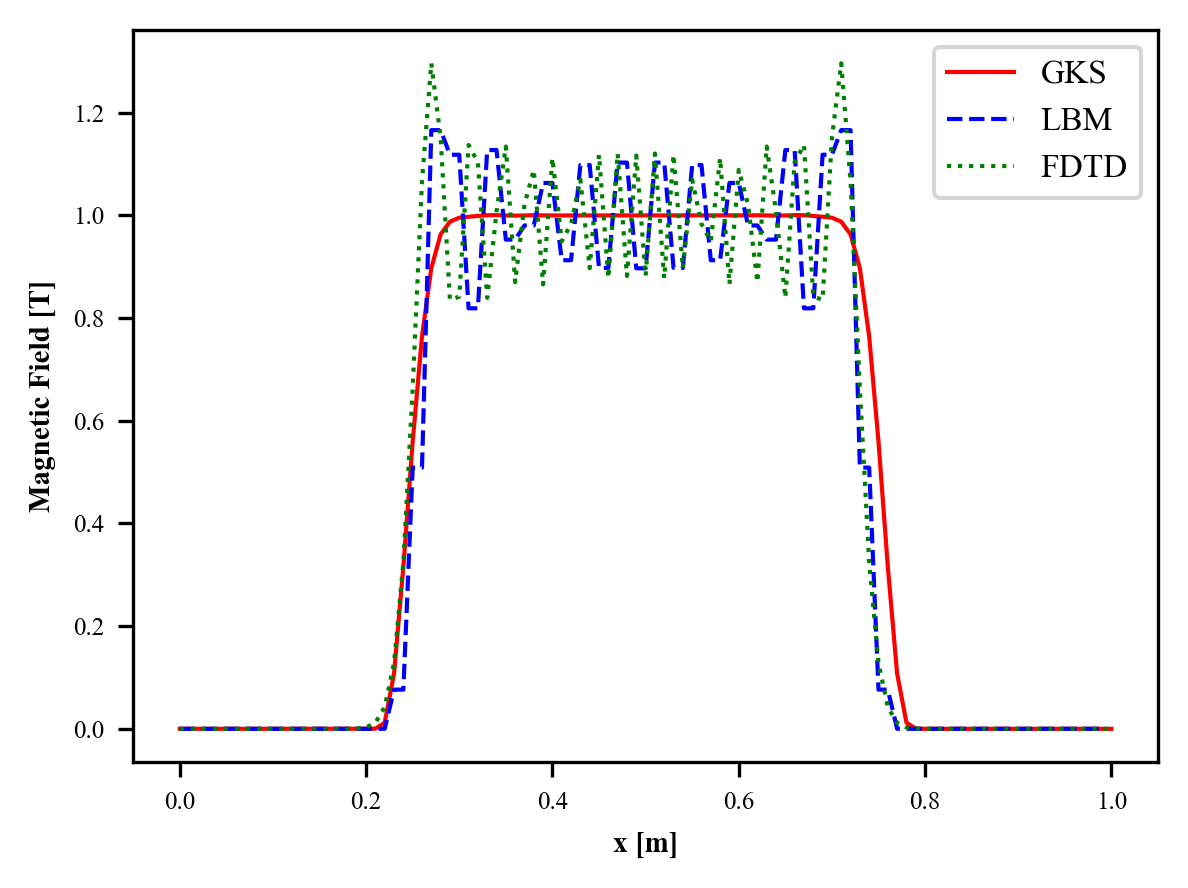} 
    \caption{\(B_z\)}
\end{subfigure}
\caption{Riemann problem: solution of the electromagnetic Riemann problem at \(t=0.25\)s comparing GKS, LBM, and FDTD methods. The GKS resolves discontinuities using a limiter. The second-order LBM and FDTD scheme produces expected oscillations near discontinuities.}
\label{fig:rp}
\end{figure}

The simulation runs until \(t = 0.25 \, \text{s}\) on a uniform grid of 100 cells with a CFL number of 0.5 for the GKS. Figure \ref{fig:rp} presents the results. The GKS employs a limiter that reduces the scheme to first-order accuracy at the discontinuity while maintaining second-order accuracy in smooth regions, thereby capturing the wave structure sharply without spurious oscillations. In contrast, the LBM maintains a global second-order accuracy by using a relaxation parameter of $\omega=2$, which renders it unstable when encountering discontinuities and leads to oscillations in these regions. As expected for a second-order central scheme, the FDTD method exhibits significant oscillations near the discontinuities. 

\subsection{Oscillating electric dipole}

This section presents a simulation of the radiation field generated by a small oscillating electric dipole antenna. The computational domain spans \([0, 1] \times [0, 1]\)m and is discretized with a \(100 \times 100\) uniform grid. The dipole is oriented along the \(z\)-axis and located at the domain center. The time-harmonic current density is prescribed as
\[
J_z = J_0 \cos\left( \frac{2\pi}{T} t \right),
\]
where \(J_0 = 1 A/m^2\) is the amplitude and \(T = 0.1\)s is the oscillation period. Non-reflecting boundary conditions are applied on all sides. The simulation is run until \(t = 3\,\text{s}\).

\begin{figure}
\centering
\begin{subfigure}[b]{0.48\textwidth}
    \includegraphics[width=0.95\textwidth]{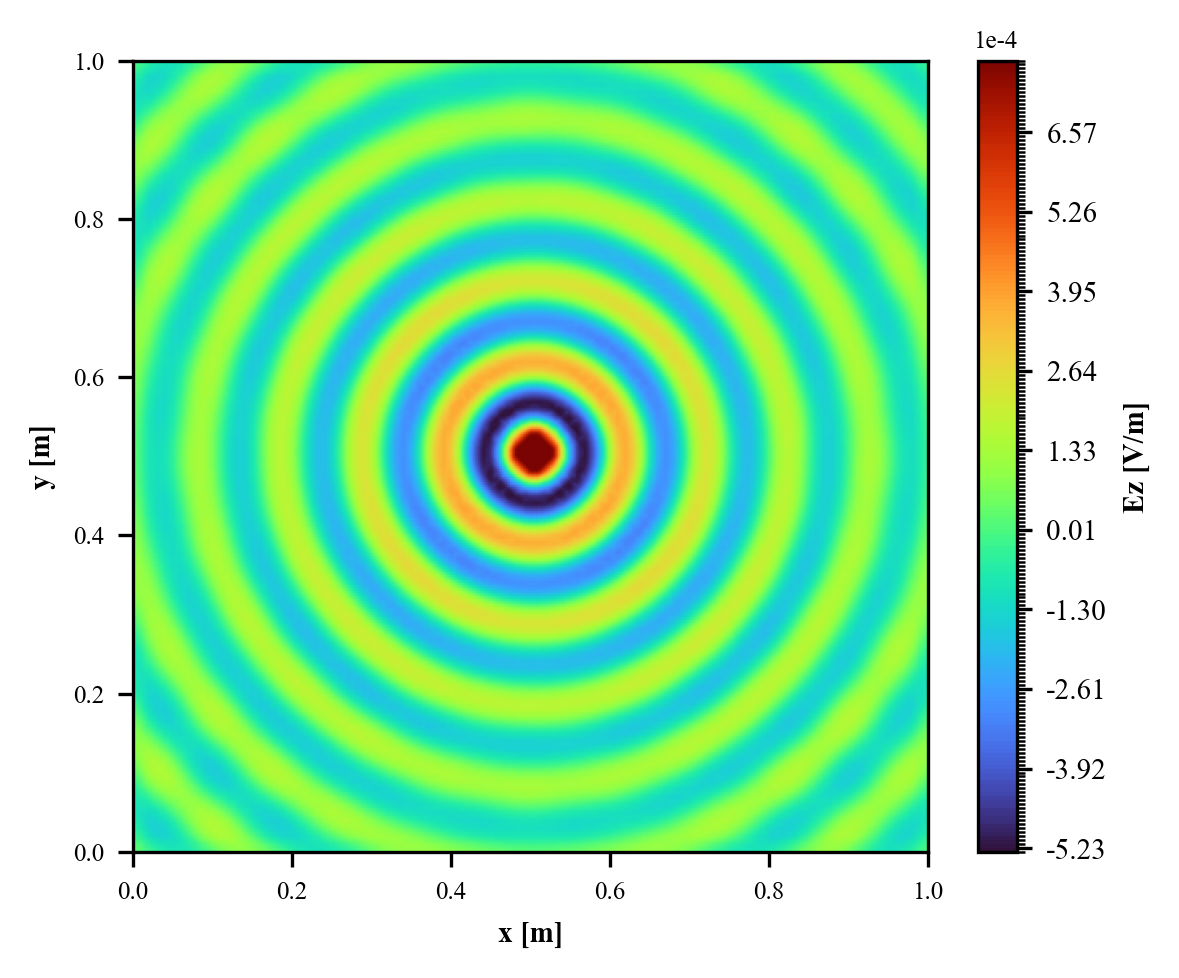}
    \caption{GKS: \(E_z\)}
\end{subfigure}
\hfill
\begin{subfigure}[b]{0.48\textwidth}
    \includegraphics[width=0.95\textwidth]{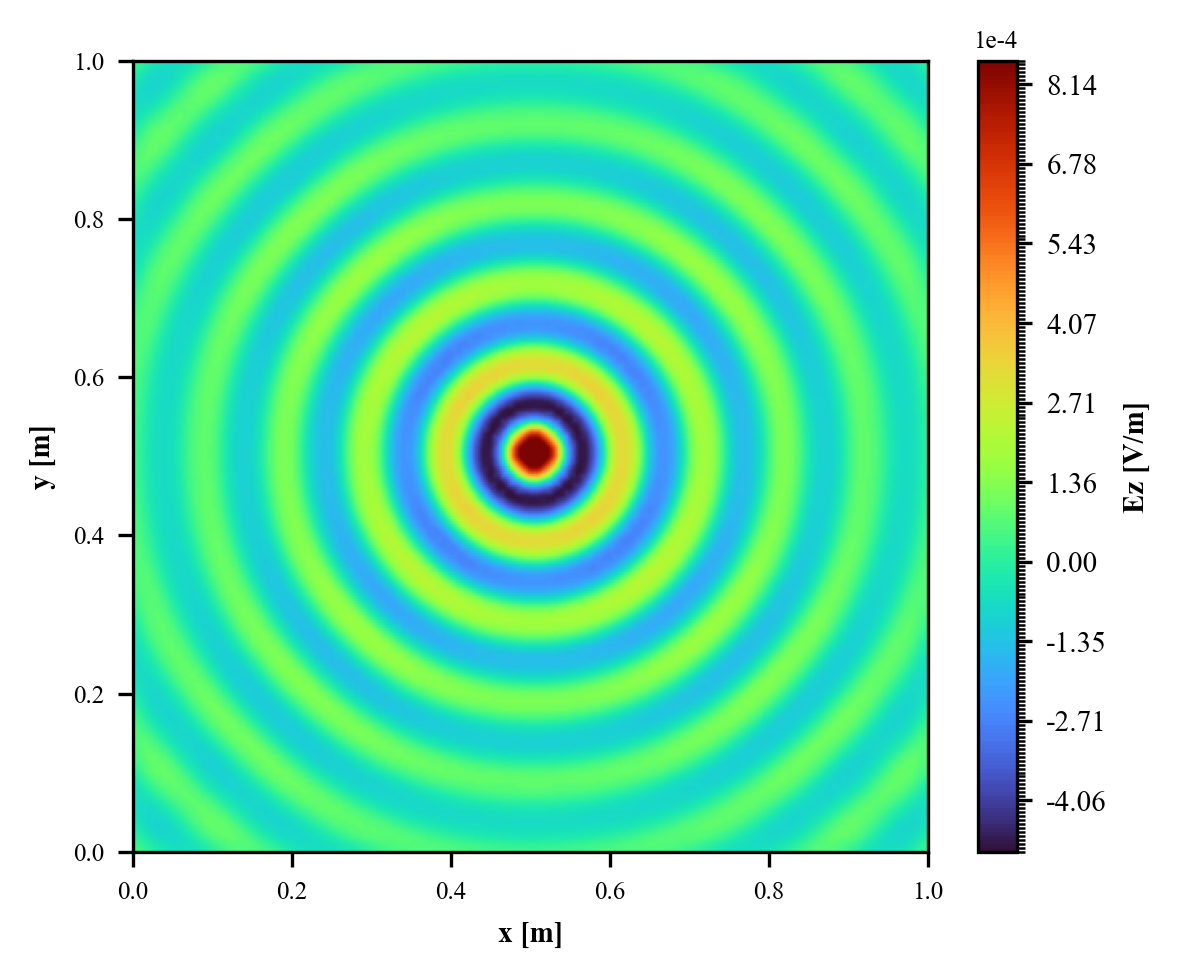}
    \caption{FVS: \(E_z\)}
\end{subfigure}
\vfill
\begin{subfigure}[b]{0.48\textwidth}
    \includegraphics[width=0.95\textwidth]{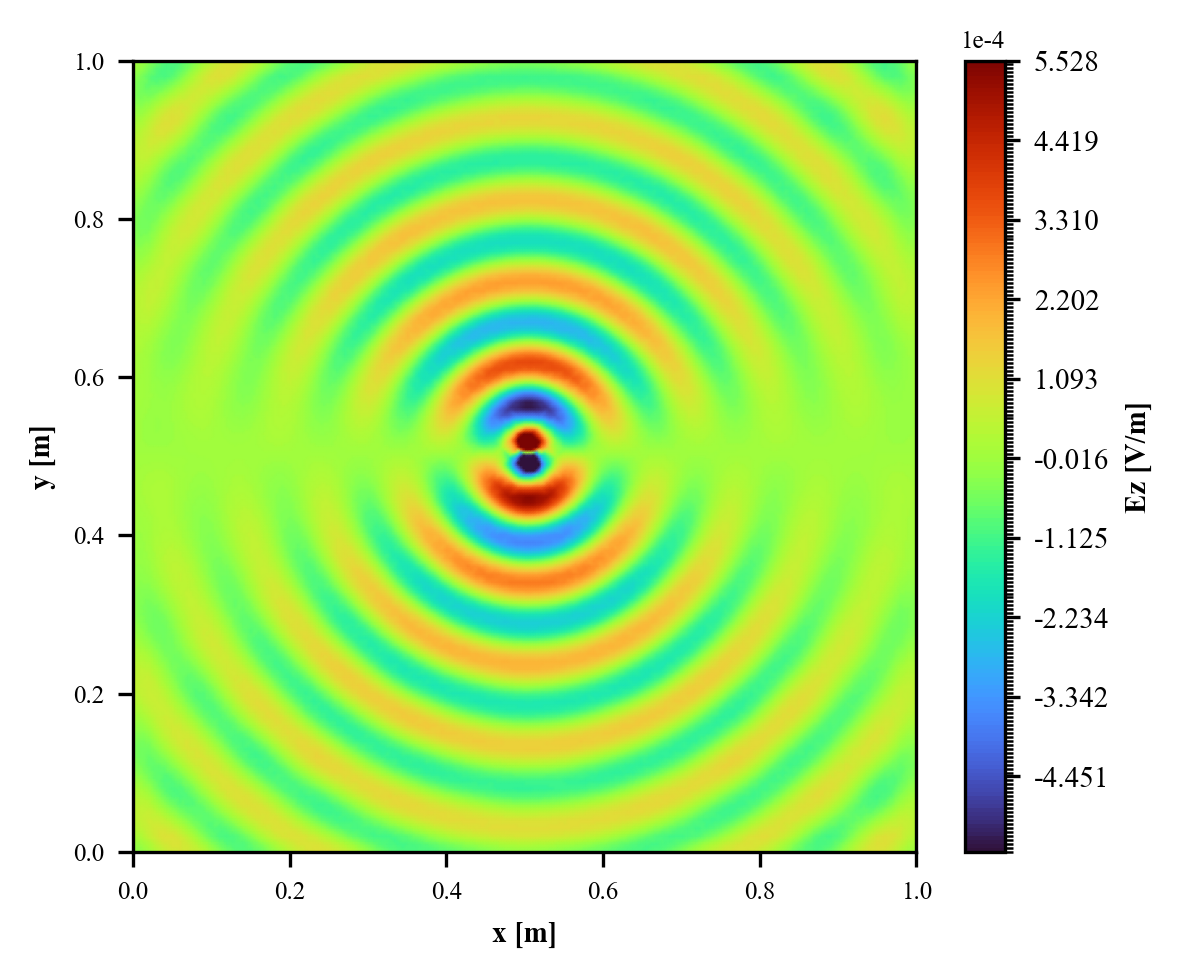}
    \caption{GKS: \(B_x\)}
\end{subfigure}
\hfill
\begin{subfigure}[b]{0.48\textwidth}
    \includegraphics[width=0.95\textwidth]{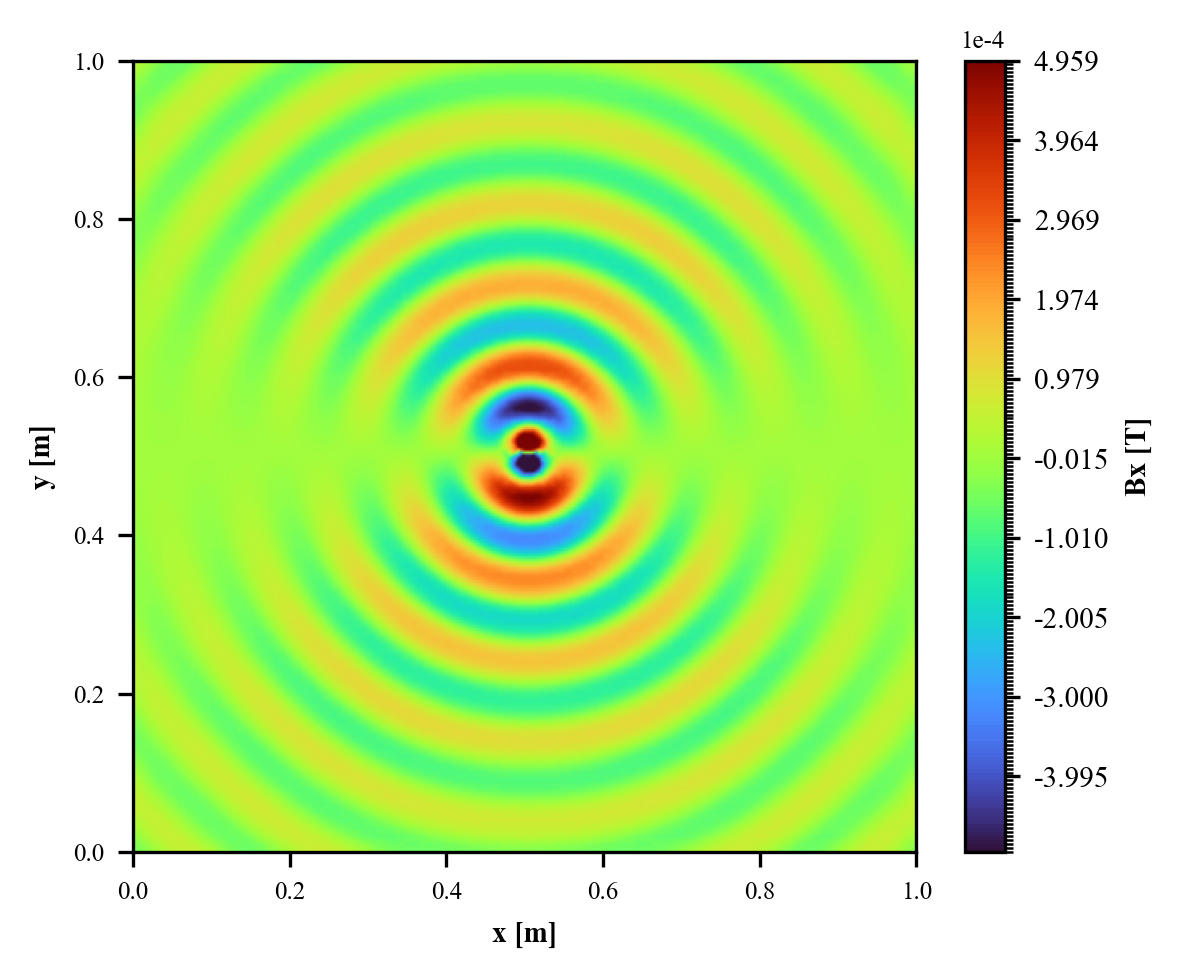}
    \caption{FVS: \(B_x\)}
\end{subfigure}
\vfill
\begin{subfigure}[b]{0.48\textwidth}
    \includegraphics[width=0.95\textwidth]{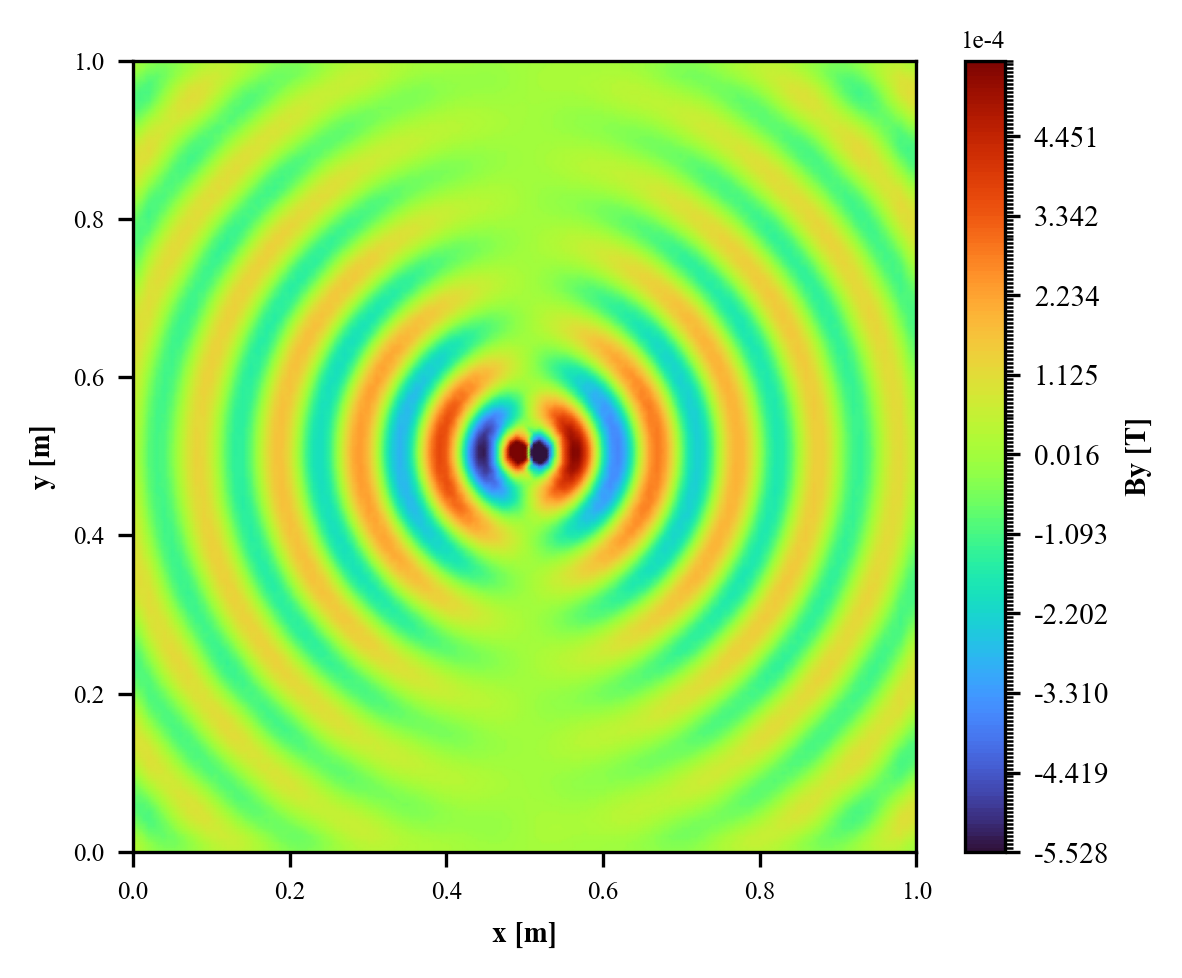}
    \caption{GKS: \(B_y\)}
\end{subfigure}
\hfill
\begin{subfigure}[b]{0.48\textwidth}
    \includegraphics[width=0.95\textwidth]{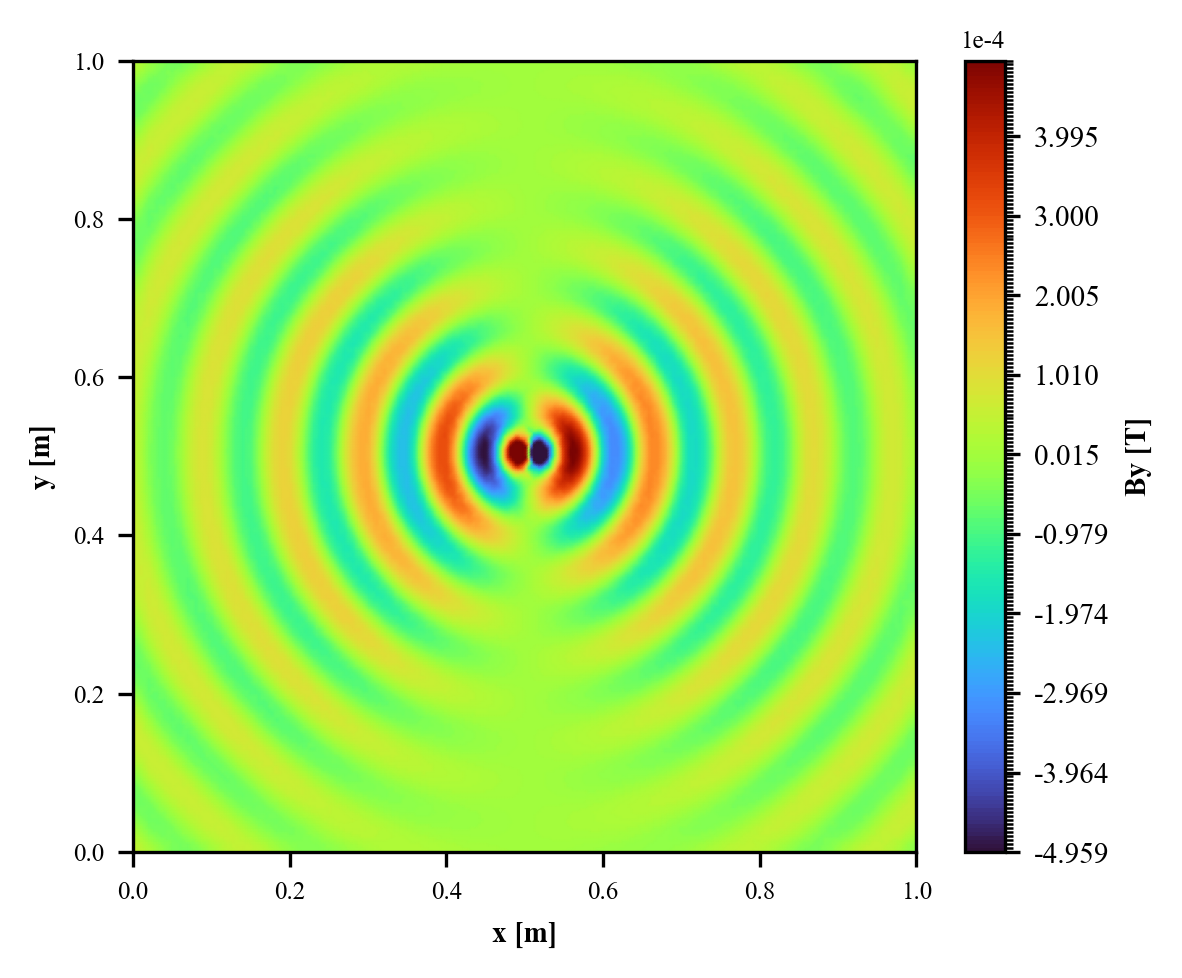}
    \caption{FVS: \(B_y\)}
\end{subfigure}
\caption{Oscillating electric dipole: Contours of the radiating fields at \(t = 3\)s for the oscillating dipole test. Left column: GKS results. Right column: FVS results. Rows from top to bottom show \(E_z\), \(B_x\), and \(B_y\).}
\label{fig:dipole1contour}
\end{figure}

Figure \ref{fig:dipole1contour} compares the field contours for GKS (left column) and FVS (right column) at the final time. The first, second, and third rows display \(E_z\), \(B_x\), and \(B_y\) contours, respectively. Both schemes capture the characteristic oscillatory wave pattern of the radiating dipole, and their results show good qualitative agreement. 

\begin{figure}
\centering
\begin{subfigure}[b]{0.48\textwidth}
    \includegraphics[width=0.95\textwidth]{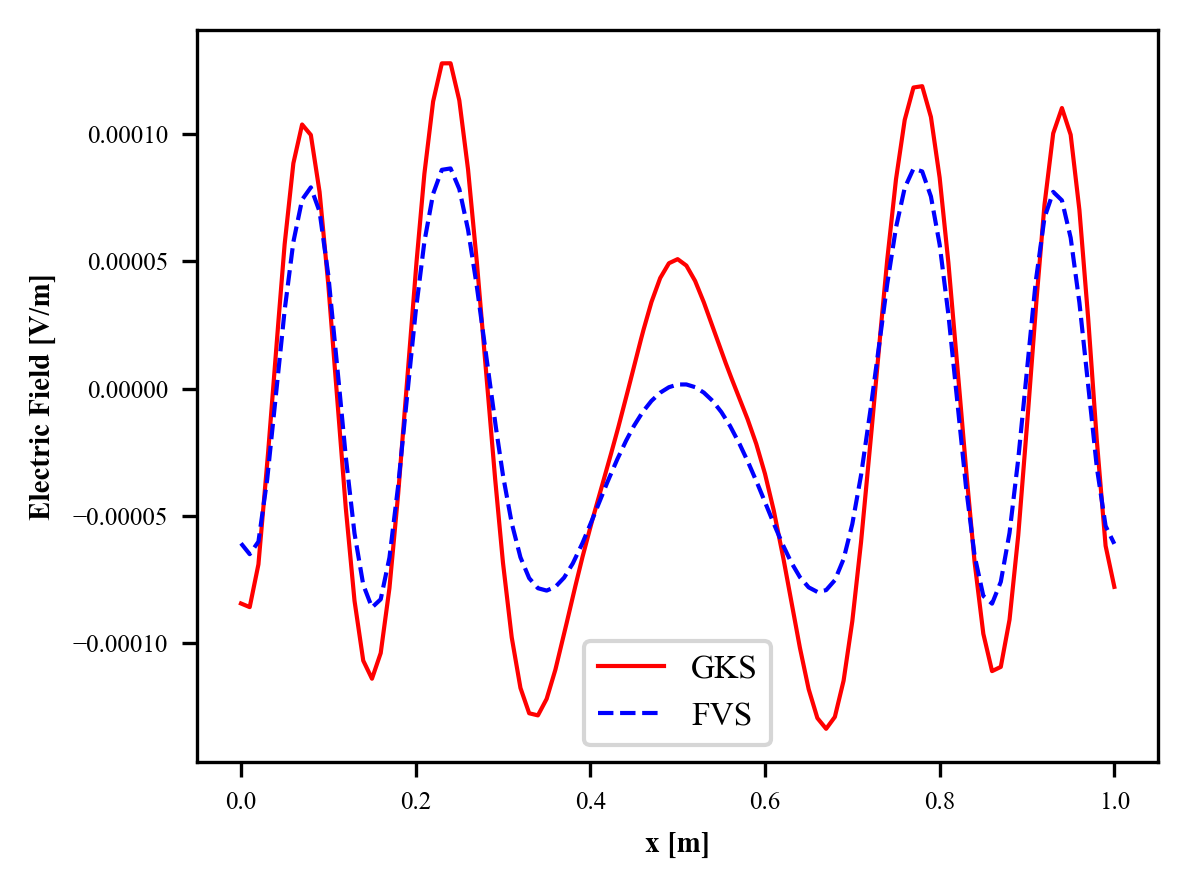}
    \caption{\(E_z\) along \(y = 0.5\) m}
\end{subfigure}
\hfill
\begin{subfigure}[b]{0.48\textwidth}
    \includegraphics[width=0.95\textwidth]{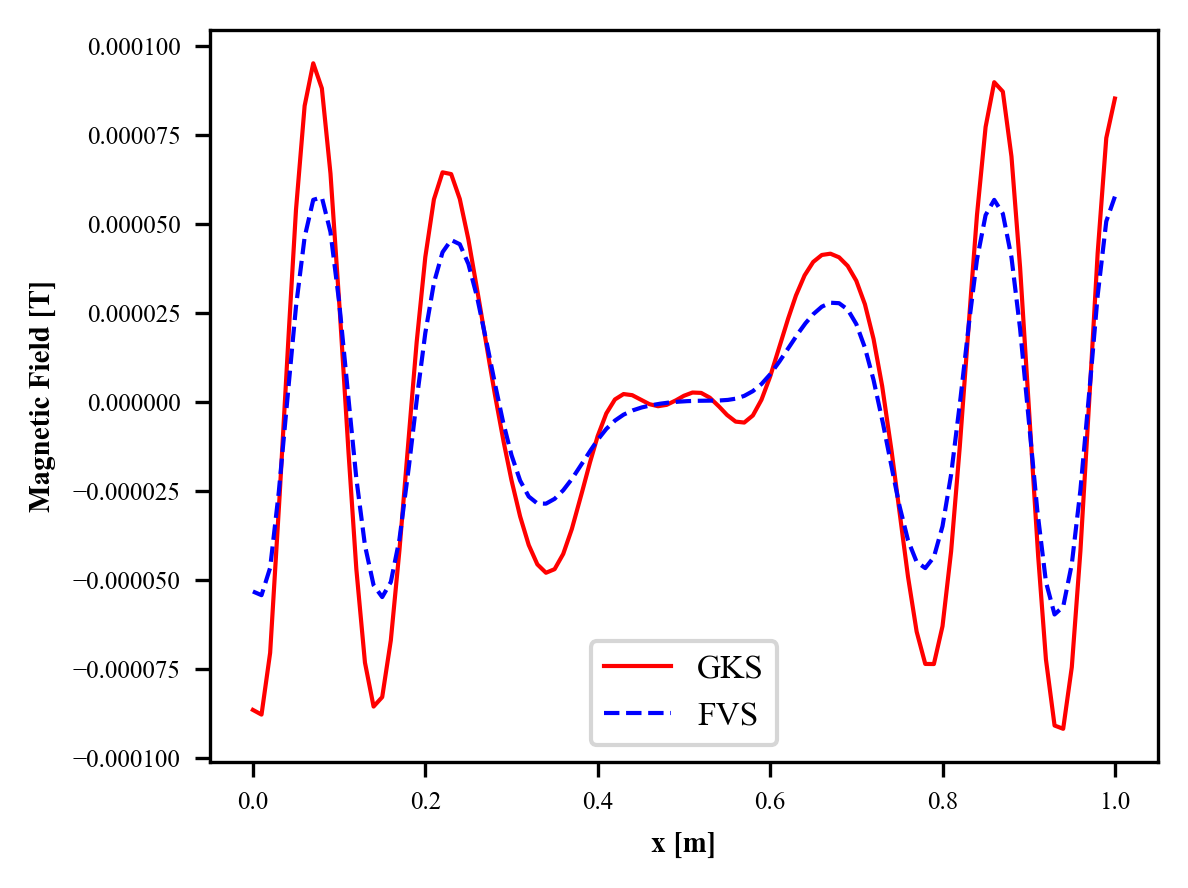}
    \caption{\(B_y\) along \(y = 0.5\) m}
\end{subfigure}
\vfill
\begin{subfigure}[b]{0.48\textwidth}
    \includegraphics[width=0.95\textwidth]{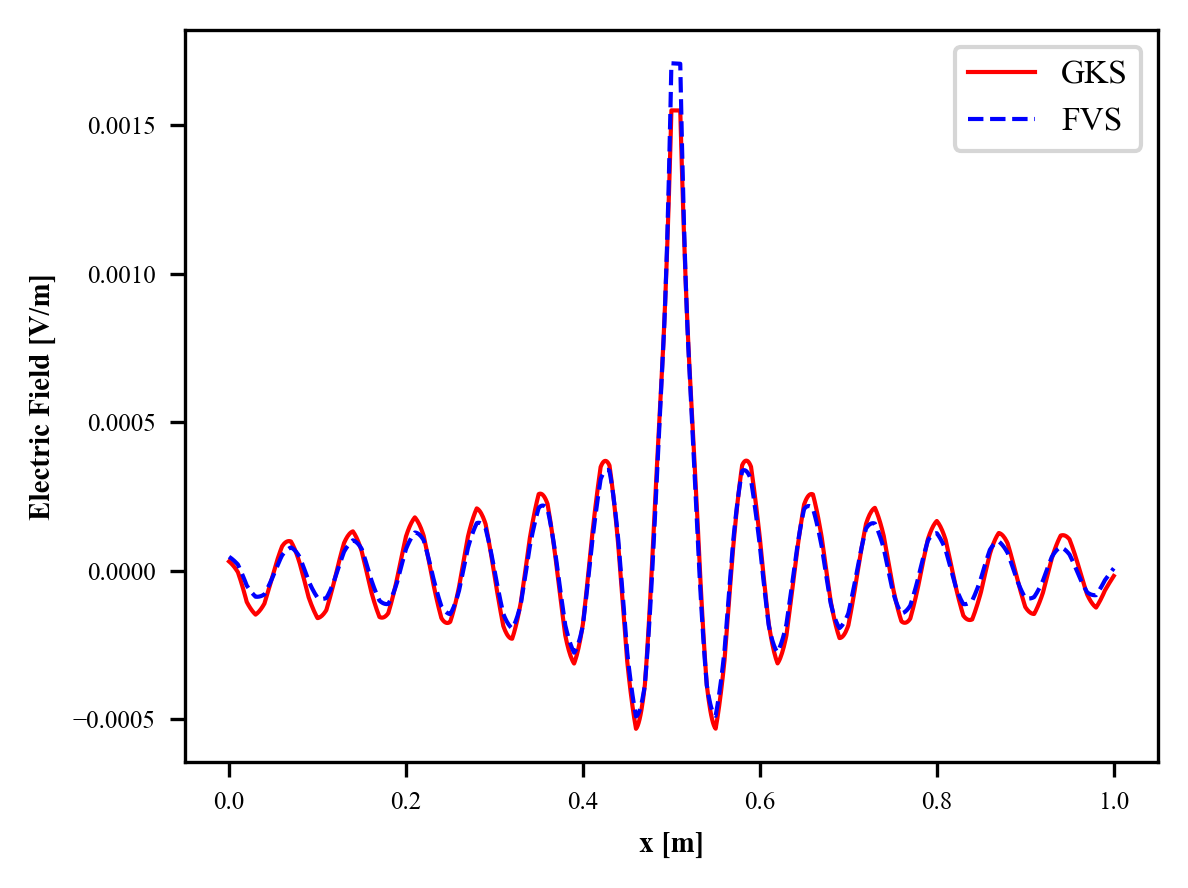}
    \caption{\(E_z\) along the diagonal}
\end{subfigure}
\hfill
\begin{subfigure}[b]{0.48\textwidth}
    \includegraphics[width=0.95\textwidth]{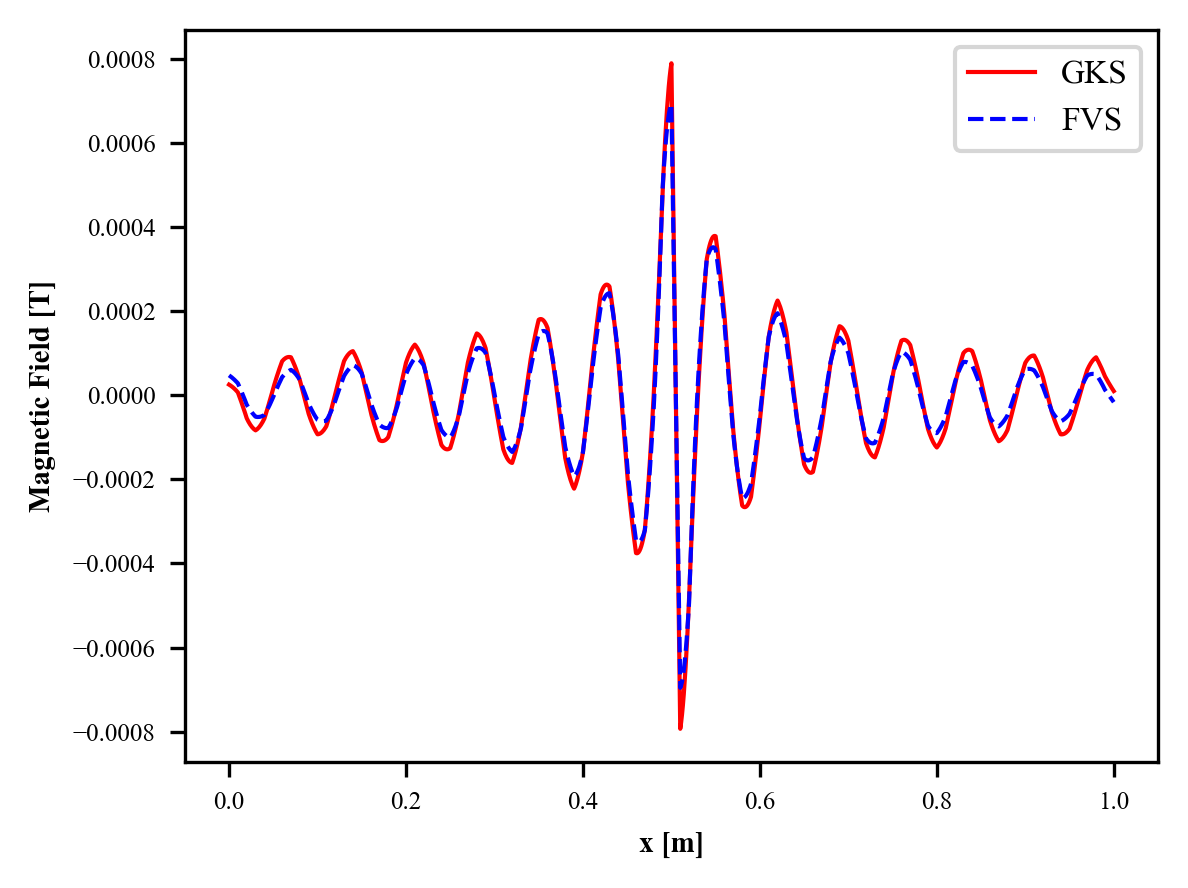}
    \caption{\(B_y\) along the diagonal}
\end{subfigure}
\caption{Oscillating electric dipole: Line profiles comparing GKS and FVS for the oscillating dipole. The FVS method shows increased numerical dissipation, leading to greater wave amplitude attenuation compared to GKS.}
\label{fig:dipole1profile}
\end{figure}

A more quantitative comparison is provided in Figure \ref{fig:dipole1profile}, which shows line profiles of \(E_z\) and \(B_y\) along the horizontal centerline (\(y = 0.5\,\text{m}\)) and the domain diagonal. The FVS results exhibit greater numerical dissipation, as evidenced by stronger damping of the wave amplitude with distance from the source. This behavior is expected, as the GKS framework naturally integrates both free-transport and collisional relaxation processes, providing a more physical representation of wave propagation than the purely upwind-based FVS approach.

We further investigate the interference pattern produced by two radiating electric dipoles. The computational domain is extended to [0, 2] $\times$ [0,1]m. The two dipoles, identical to the one described previously, are positioned at (0.7, 0.5)m and (1.3, 0.5)m, respectively. The simulation is run until \(t = 2\,\text{s}\).

\begin{figure}
    \centering
    \includegraphics[width=0.9\linewidth]{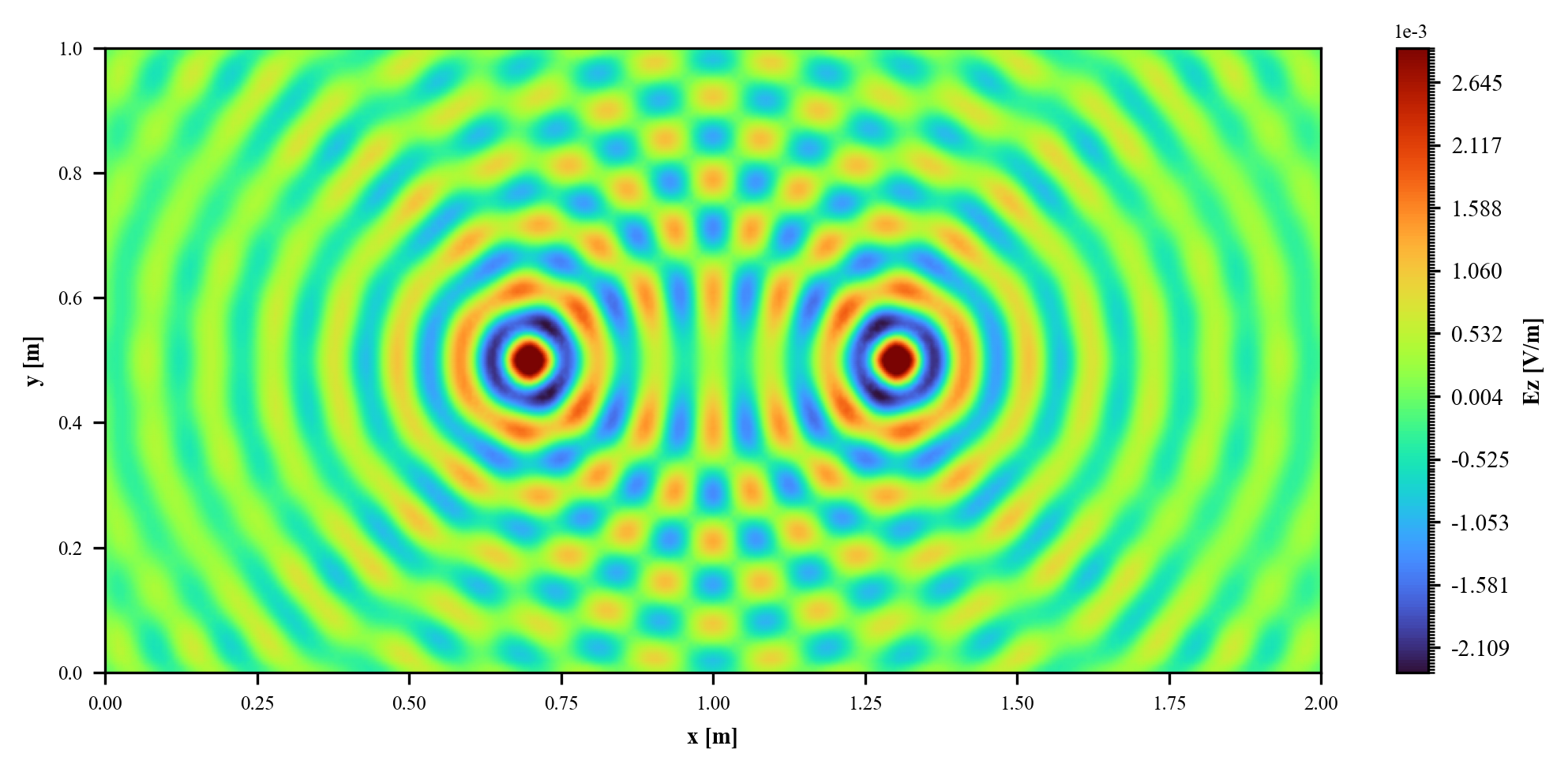}
    \caption{Oscillating electric dipole: $E_z$ contour by GKS.}
    \label{fig:dipole2Ez}
\end{figure}

\begin{figure}
    \centering
    \includegraphics[width=0.9\linewidth]{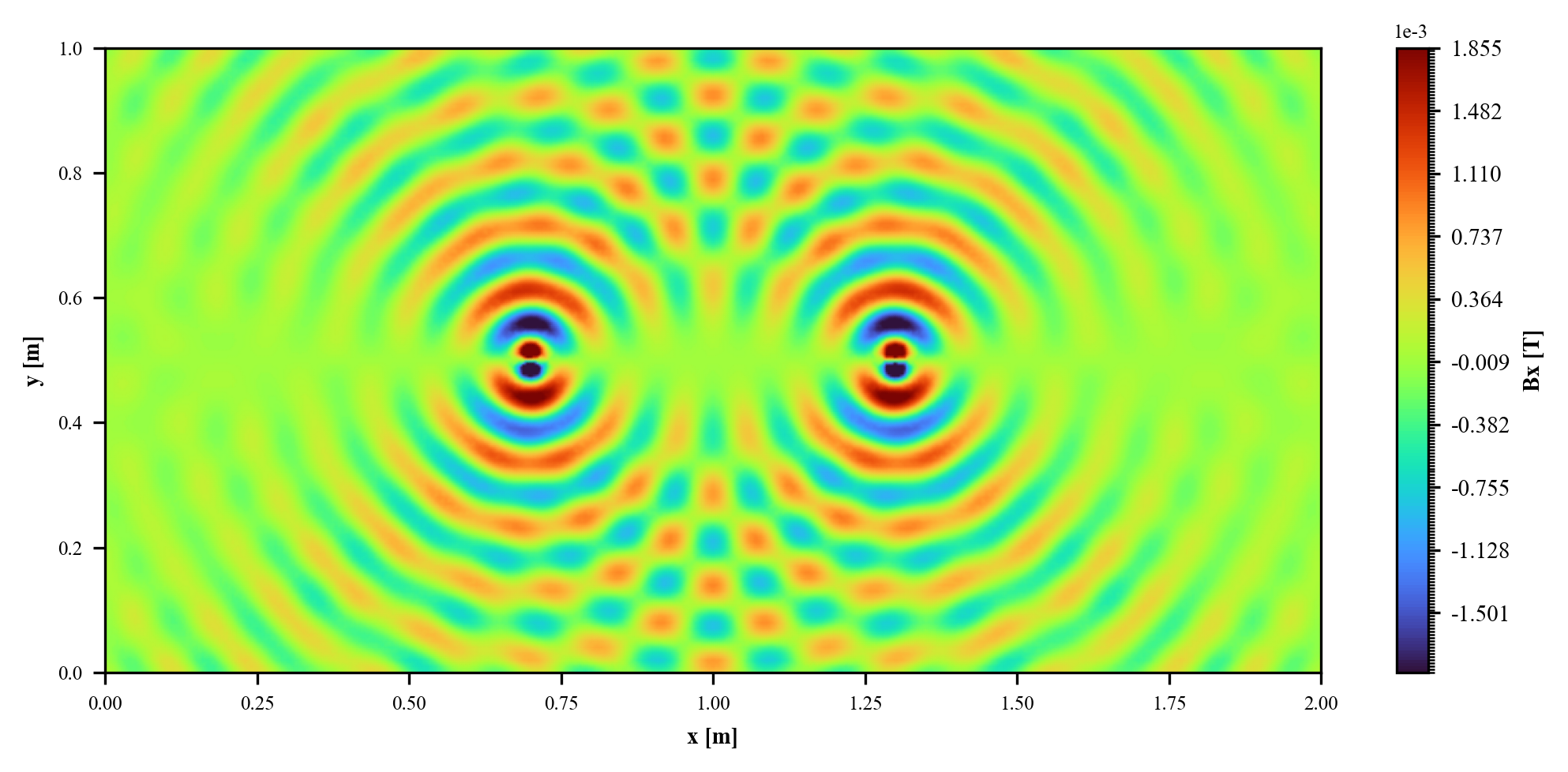}
    \caption{Oscillating electric dipole: $B_x$ contour by GKS.}
    \label{fig:dipole2Bx}
\end{figure}

\begin{figure}
    \centering
    \includegraphics[width=0.9\linewidth]{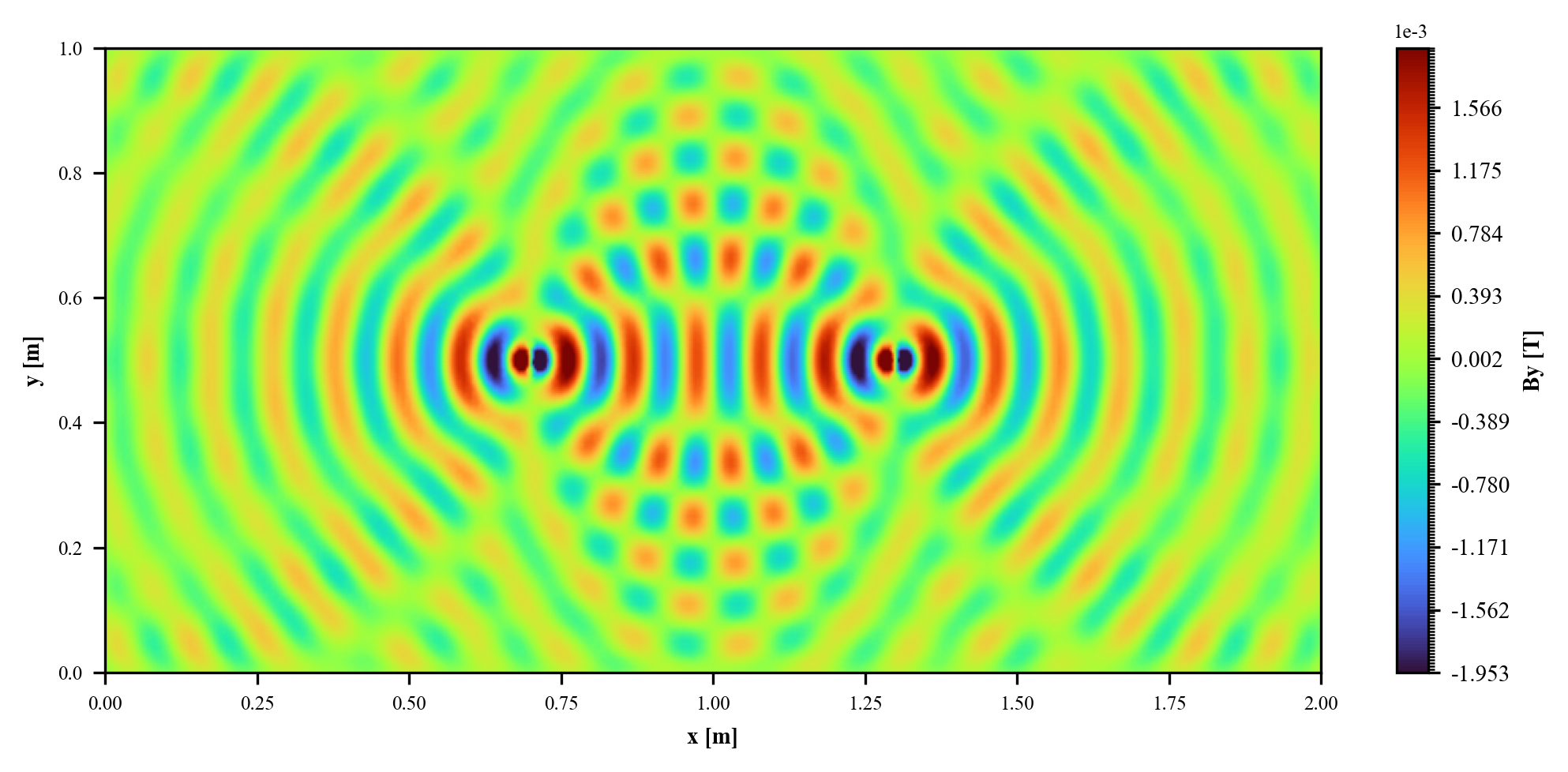}
    \caption{Oscillating electric dipole: $B_y$ contour by GKS.}
    \label{fig:dipole2By}
\end{figure}

\begin{figure}
\centering
\begin{subfigure}[b]{0.48\textwidth}
    \includegraphics[width=0.95\textwidth]{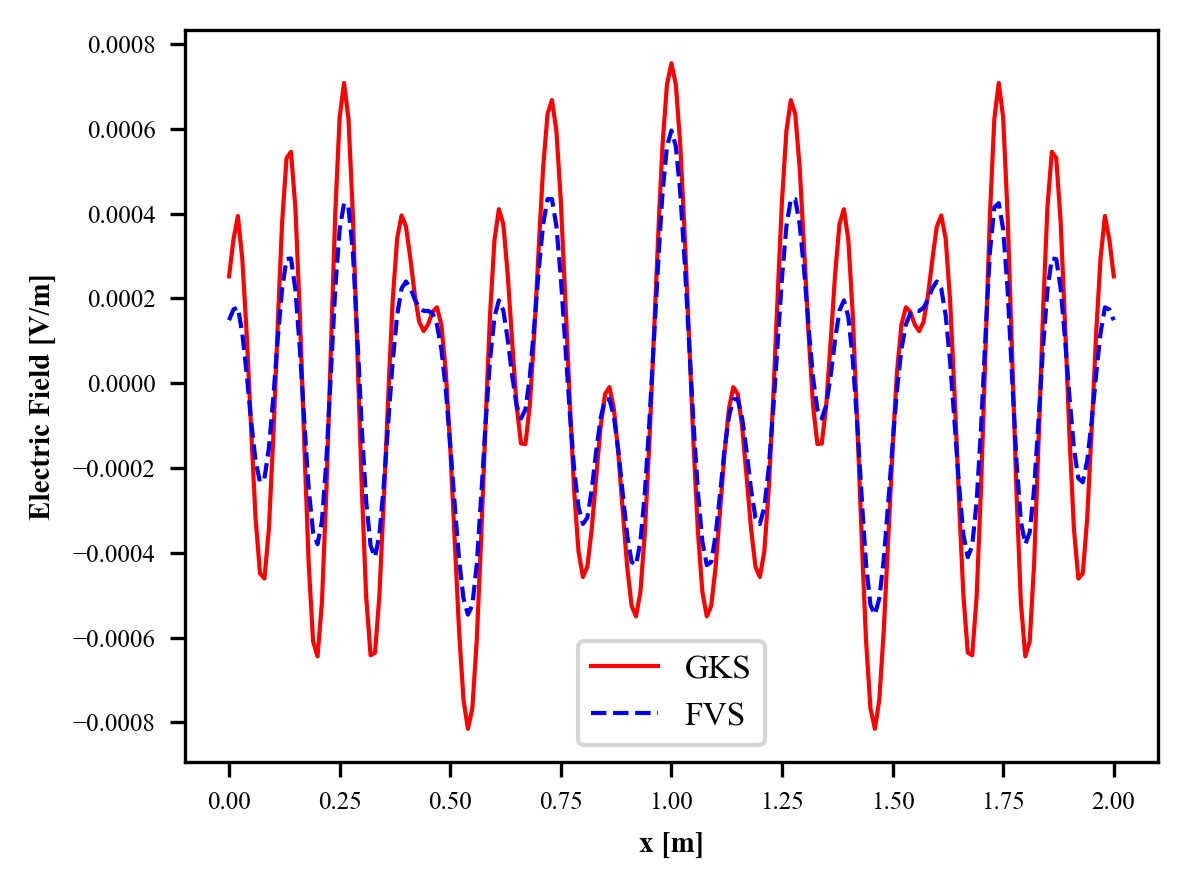}
    \caption{\(E_z\) along \(y = 0.5\) m}
\end{subfigure}
\hfill
\begin{subfigure}[b]{0.48\textwidth}
    \includegraphics[width=0.95\textwidth]{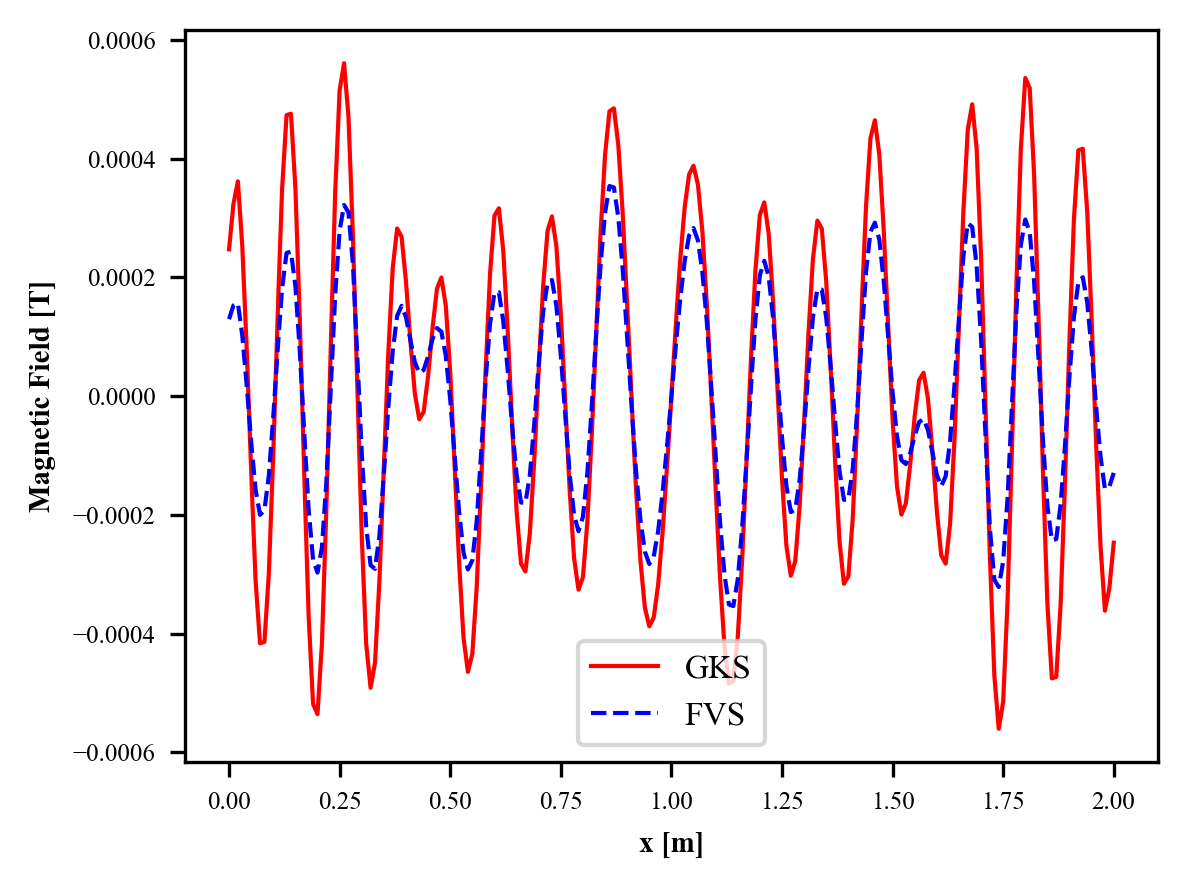}
    \caption{\(B_y\) along \(y = 0.5\) m}
\end{subfigure}
\vfill
\begin{subfigure}[b]{0.48\textwidth}
    \includegraphics[width=0.95\textwidth]{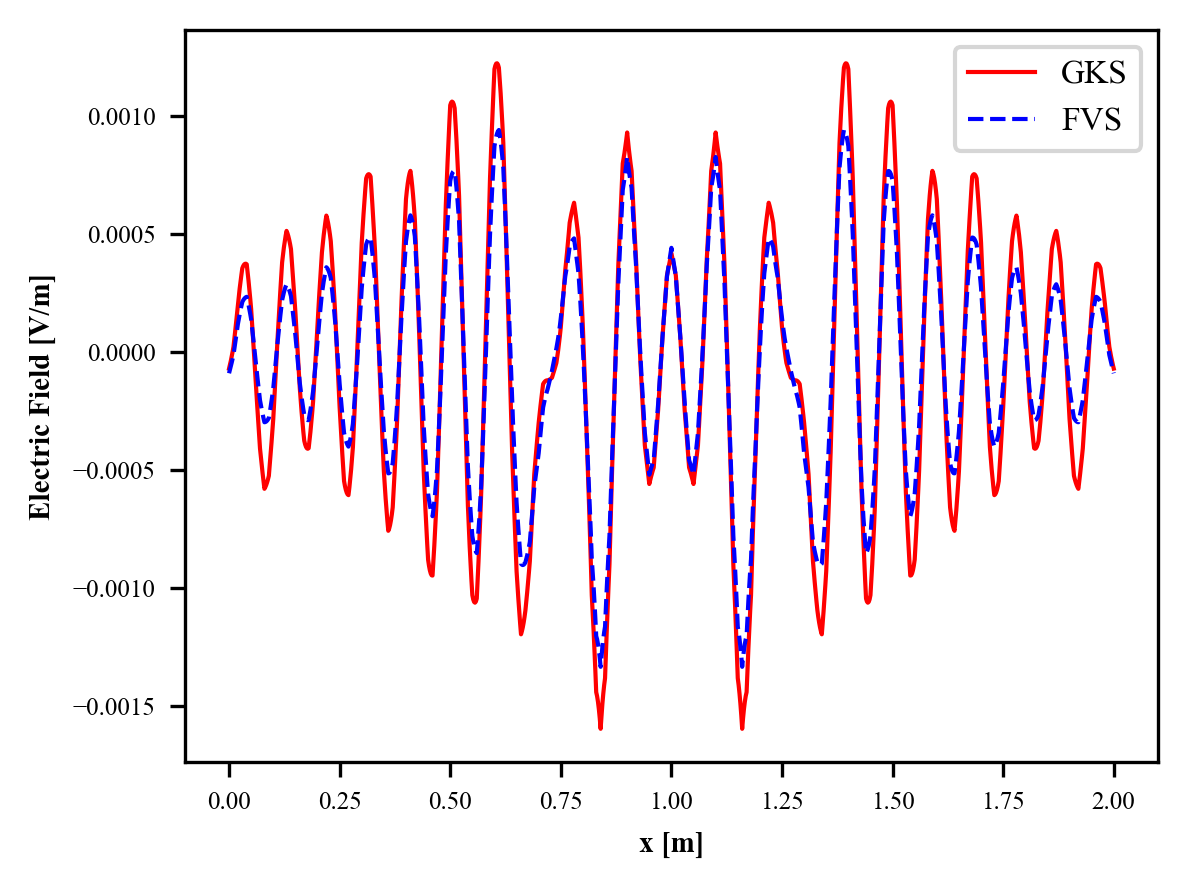}
    \caption{\(E_z\) along  the diagonal}
\end{subfigure}
\hfill
\begin{subfigure}[b]{0.48\textwidth}
    \includegraphics[width=0.95\textwidth]{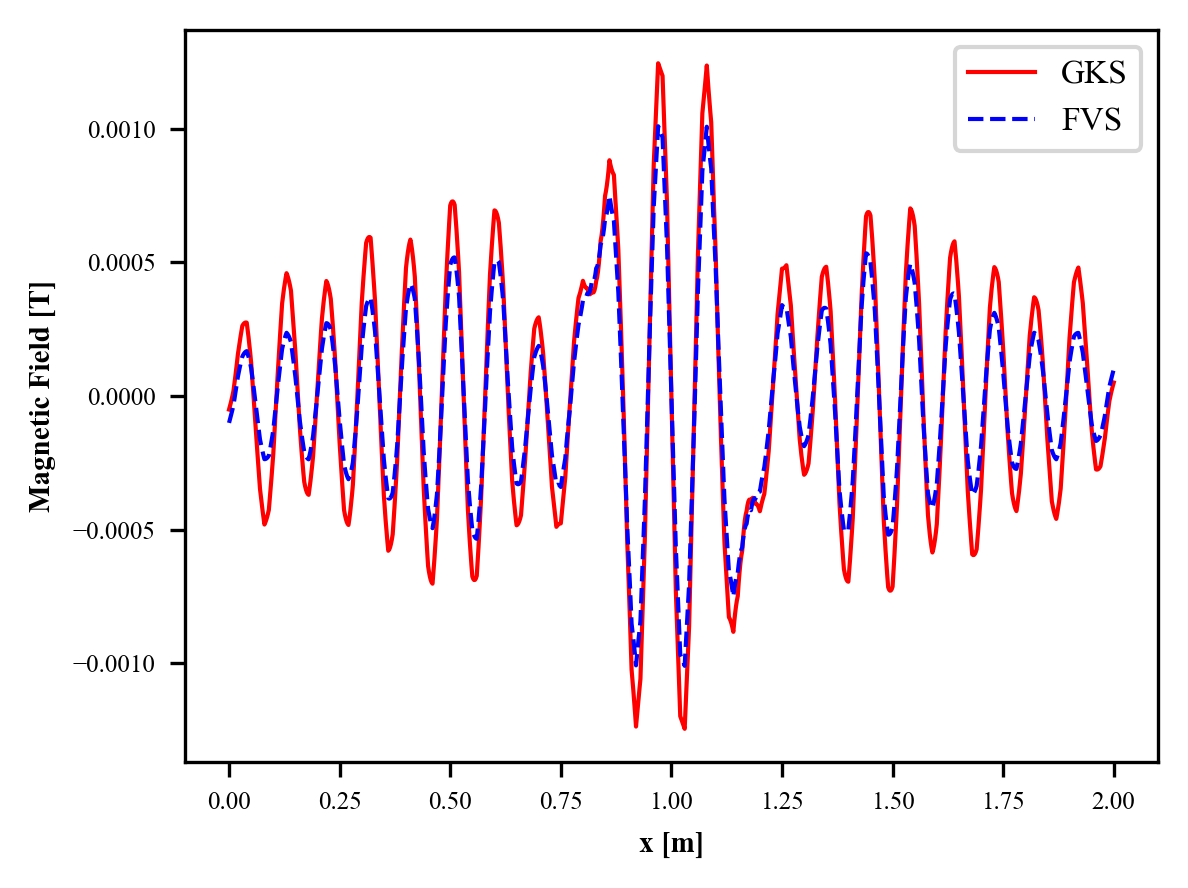}
    \caption{\(B_y\) along  the diagonal}
\end{subfigure}
\caption{Oscillating electric dipole: Line profiles comparing GKS and FVS for the two oscillating dipoles. The FVS method shows increased numerical dissipation, leading to greater wave amplitude attenuation compared to GKS.}
\label{fig:dipole2profile}
\end{figure}

The resulting radiation fields computed with the GKS are shown in Figures \ref{fig:dipole2Ez}, \ref{fig:dipole2Bx}, and \ref{fig:dipole2By}, which display contours of \(E_z\), \(B_x\), and \(B_y\), respectively. The expected interference pattern between the two sources is accurately captured. Line profiles along the horizontal centerline and the domain diagonal are presented in Figure \ref{fig:dipole2profile}. The solutions from both the GKS and FVS methods are closely, and the comparative dissipation trends observed in the single-dipole case are consistent here.

\subsection{Signal propagation around an aircraft}

This section demonstrates the GKS's capability to simulate electromagnetic wave propagation around a complex geometry using an unstructured mesh. The test case models a small antenna, represented as a surface electric dipole, positioned on the aircraft fuselage. The time-harmonic current density of the dipole is
\[
J_z = J_0 \cos\left( \frac{2\pi}{T} t \right),
\]
where \(J_0 = 1 A/m^2\) and \(T = 0.1\)s are the amplitude and period, respectively. Simulations are performed with the antenna placed at several different locations to compare the resulting radiation patterns.

The aircraft surface is modeled as a perfect electric conductor (PEC). The general interface conditions for electromagnetic fields are given by
\begin{equation}
\left\{
\begin{aligned}
    &\boldsymbol{n}\times(\boldsymbol{E}_1 - \boldsymbol{E}_2) = 0,\\
    &\boldsymbol{n}\times(\boldsymbol{H}_1 - \boldsymbol{H}_2) = \boldsymbol{J}_s,\\
    &\boldsymbol{n}\cdot(\boldsymbol{B}_1 - \boldsymbol{B}_2) = 0,\\
    &\boldsymbol{n}\cdot(\boldsymbol{D}_1 - \boldsymbol{D}_2) = \rho_s,
\end{aligned}
\right.
\end{equation}
where \(\boldsymbol{n}\) is the unit normal to the surface, \(\boldsymbol{J}_s\) is the surface current density, and \(\rho_s\) is the surface charge density. The subscripts \(1\) and \(2\) denote values on either side of the interface. On a PEC surface, the tangential electric field and the normal magnetic flux density vanish. The unknown surface currents and charges are treated as jumps across the boundary \cite{shang1995scattered}. Consequently, the gradient of these jumps in the normal direction must be zero, leading to the Neumann conditions:
\begin{equation}
\left\{
\begin{aligned}
    &\boldsymbol{n}\cdot\nabla\left[\boldsymbol{n}\times(\boldsymbol{H}_1 - \boldsymbol{H}_2)\right] = \boldsymbol{n}\cdot\nabla\boldsymbol{J}_s = 0,\\
    &\boldsymbol{n}\cdot\nabla\left[\boldsymbol{n}\cdot(\boldsymbol{D}_1 - \boldsymbol{D}_2)\right] = \boldsymbol{n}\cdot\nabla (\rho_s/\epsilon) = 0.
\end{aligned}
\right.
\end{equation}
These conditions are applied to the normal component of the electric field and the tangential component of the magnetic field at the aircraft surface. A non-reflecting boundary condition is applied at the outer domain limits. The simulation uses a CFL number of 0.5.

\begin{figure}
\centering
\begin{subfigure}[b]{0.7\textwidth}
    \includegraphics[width=0.9\textwidth]{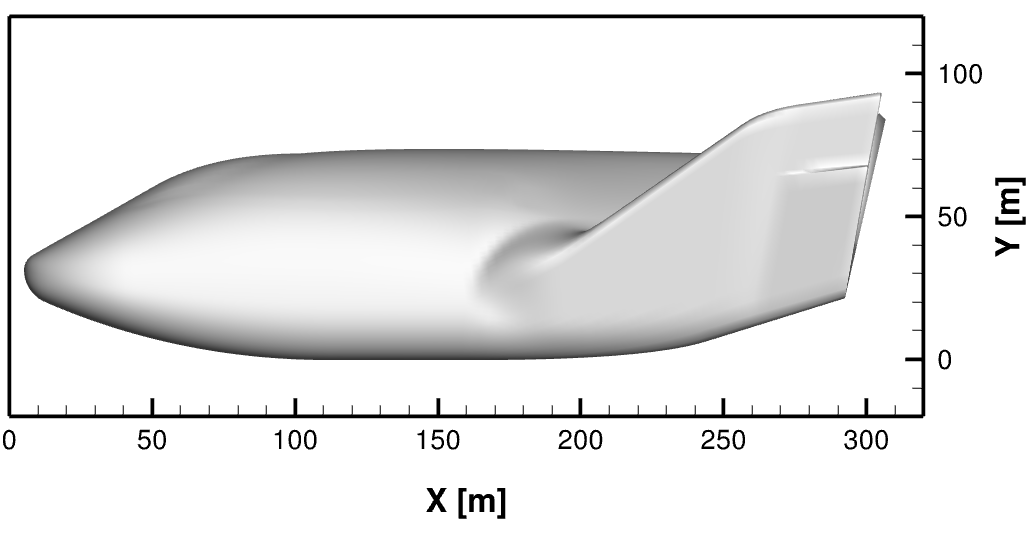} 
\end{subfigure}
\vfill
\begin{subfigure}[b]{0.45\textwidth}
    \includegraphics[width=0.9\textwidth]{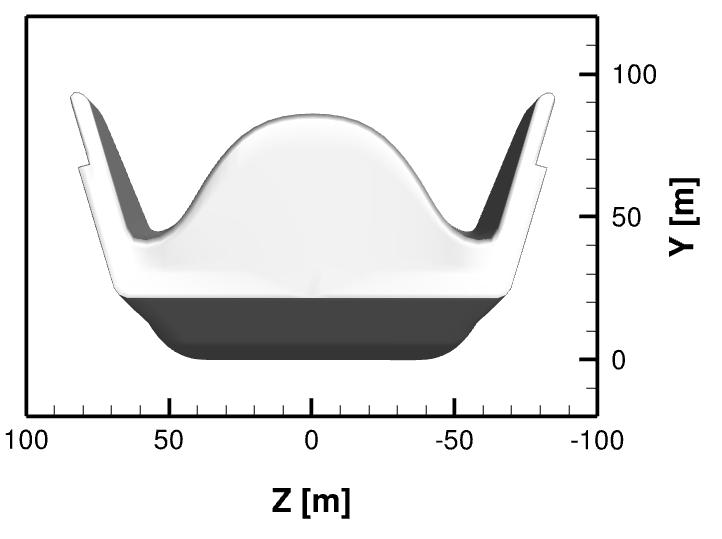} 
\end{subfigure}
\hfill
\begin{subfigure}[b]{0.45\textwidth}
    \includegraphics[width=0.98\textwidth]{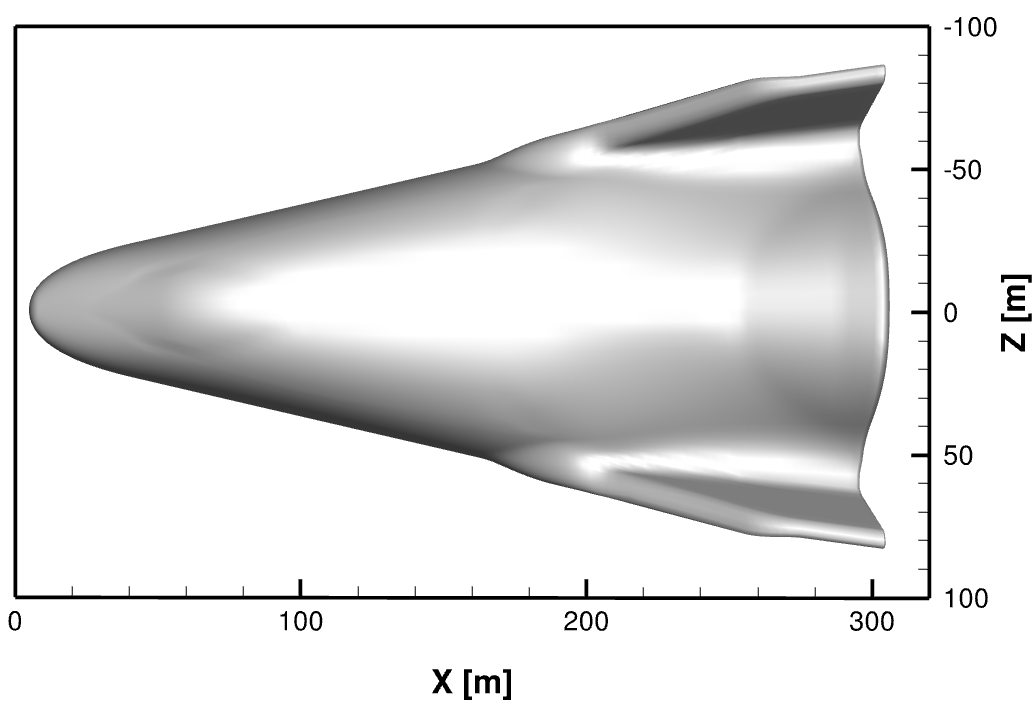} 
\end{subfigure}
\caption{Sketch of the X38-like aircraft.}
\label{fig:x38geo}
\end{figure}

\begin{figure}
\centering
\begin{subfigure}[b]{0.44\textwidth}
    \includegraphics[width=0.98\textwidth]{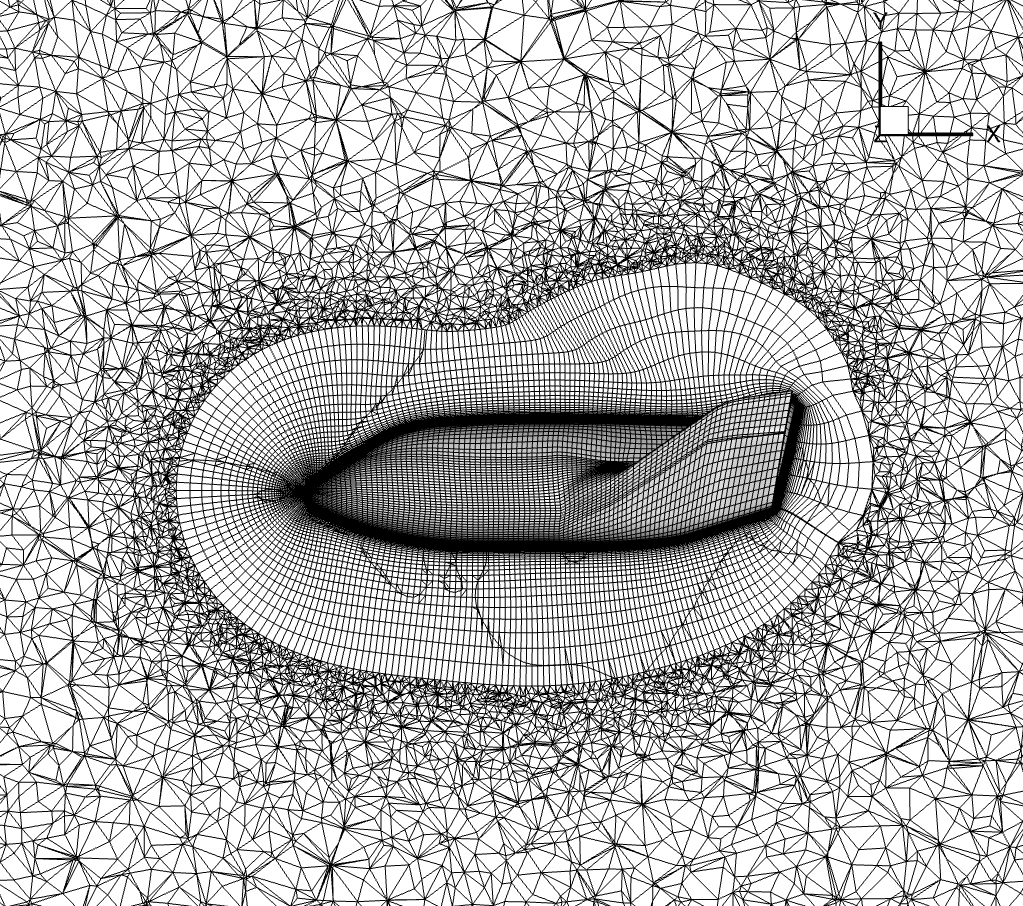} 
\end{subfigure}
\hfill
\begin{subfigure}[b]{0.48\textwidth}
    \includegraphics[width=0.98\textwidth]{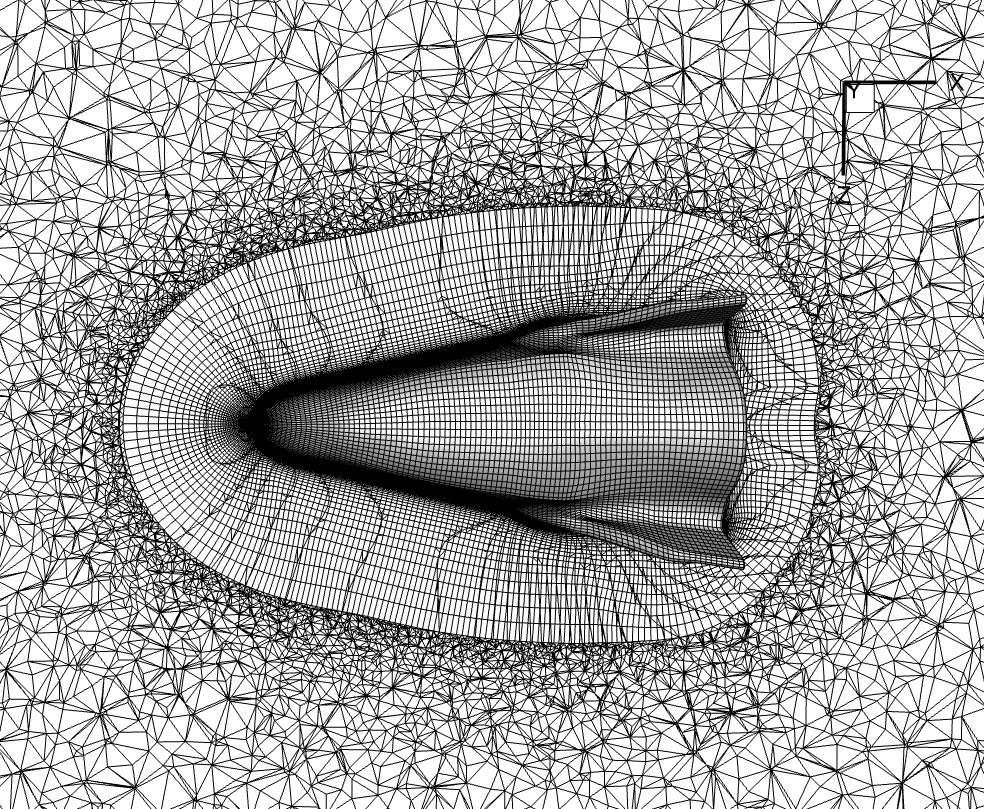} 
\end{subfigure}
\caption{Meshes around the aircraft.}
\label{fig:x38mesh}
\end{figure}

The aircraft geometry and dimensions are depicted in Figure \ref{fig:x38geo}. The computational mesh used for the simulation is shown in Figure \ref{fig:x38mesh}. The domain extends over $[-400, 1200] \, \text{m} \times [-295, 695] \, \text{m} \times [-400, 400] \, \text{m}$ and is discretized with a total of 1,685,502 cells. The mesh is refined near the aircraft surface. The simulation runs until $t=100s$.

\begin{figure}
\centering
\begin{subfigure}[b]{0.48\textwidth}
    \includegraphics[width=0.95\textwidth]{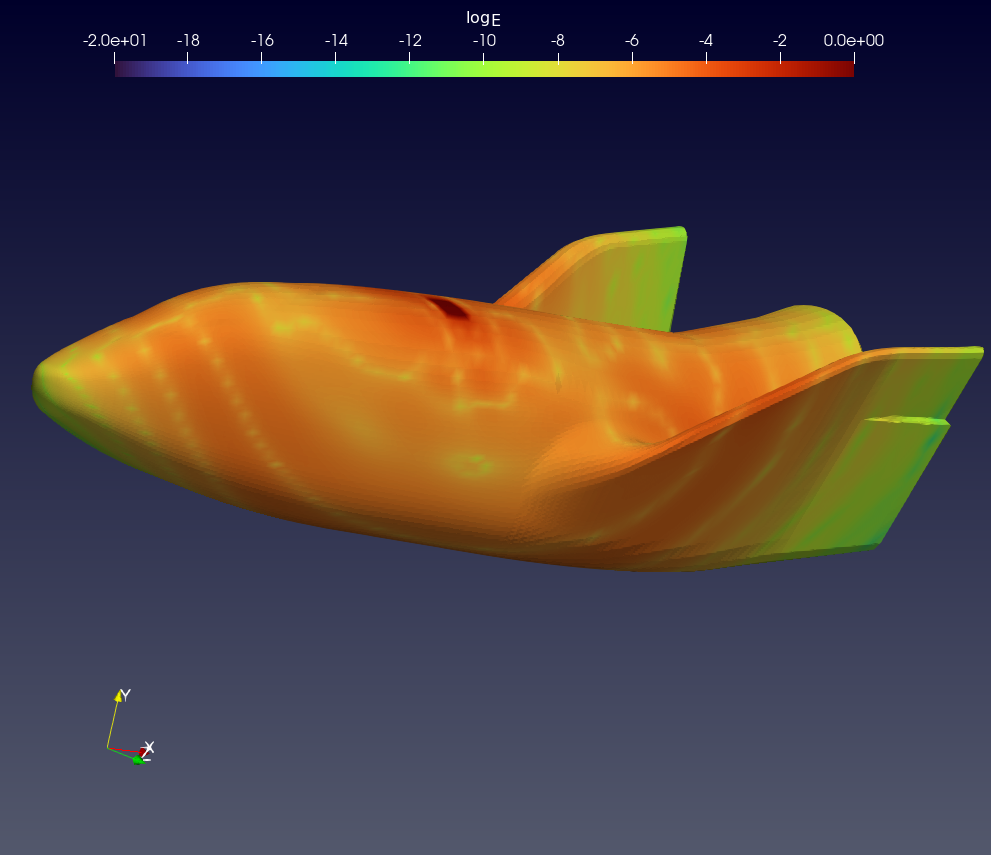}
    \caption{}
\end{subfigure}
\hfill
\begin{subfigure}[b]{0.48\textwidth}
    \includegraphics[width=0.95\textwidth]{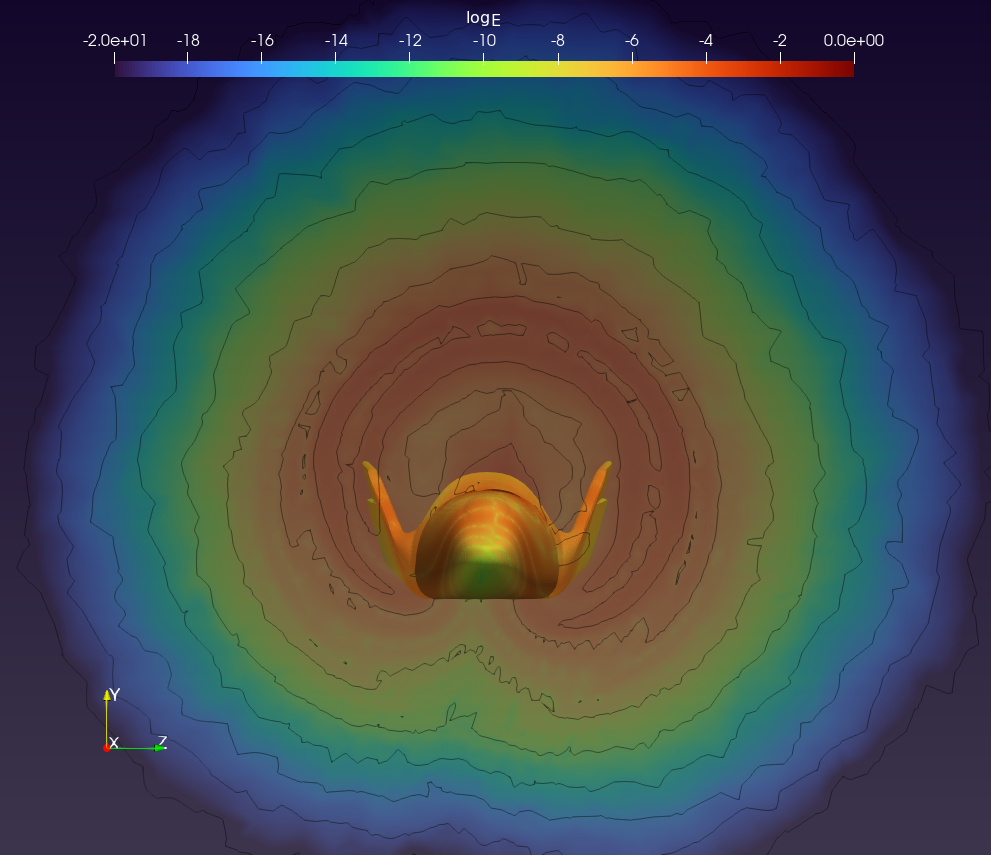}
    \caption{}
\end{subfigure}
\vfill
\begin{subfigure}[b]{0.48\textwidth}
    \includegraphics[width=0.95\textwidth]{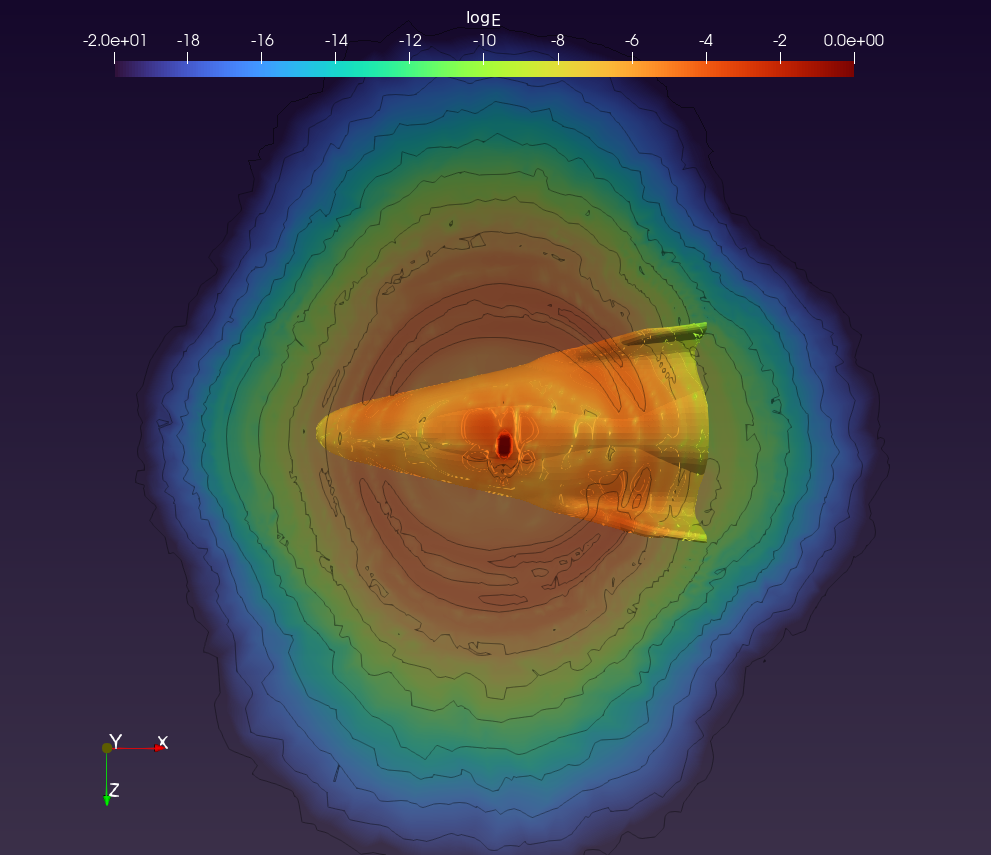}
    \caption{}
\end{subfigure}
\hfill
\begin{subfigure}[b]{0.48\textwidth}
    \includegraphics[width=0.95\textwidth]{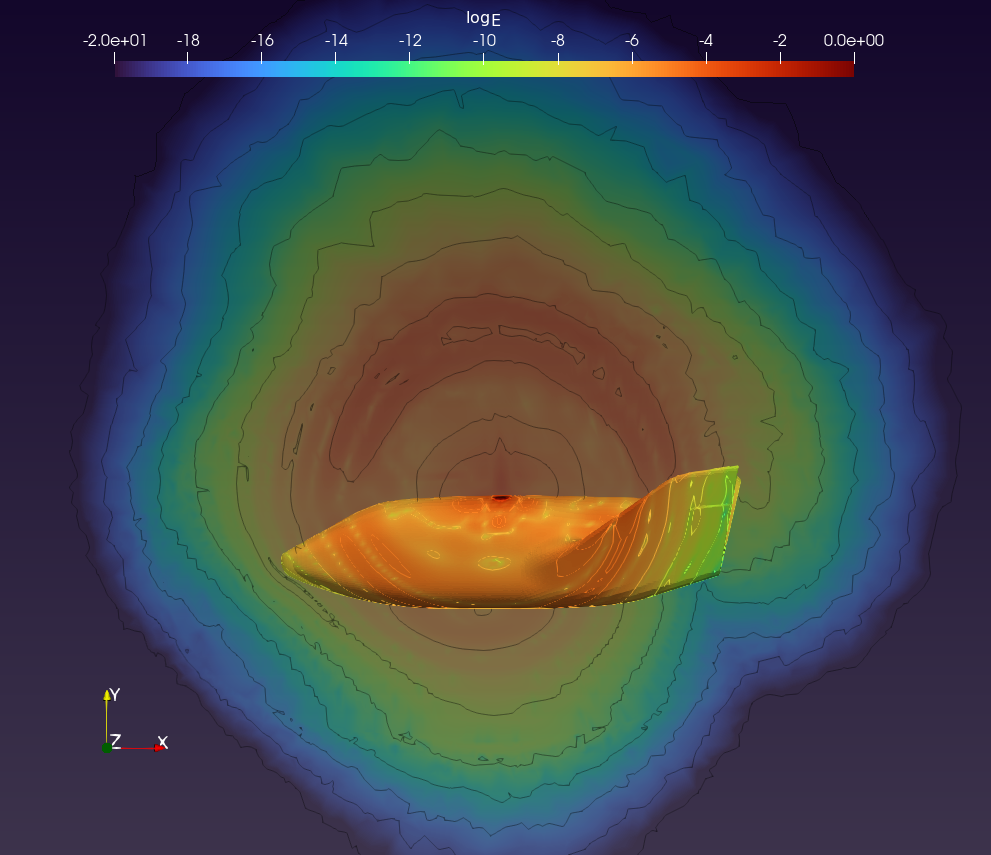}
    \caption{}
\end{subfigure}
\caption{Signal Propagation Around an Aircraft: The electric field magnitude on the surface of  (a) aircraft, (b) y-z slice, (c) x-z slice, (d) x-y slice.}
\label{fig:x381contour}
\end{figure}

\begin{figure}
\centering
\begin{subfigure}[b]{0.48\textwidth}
    \includegraphics[width=0.95\textwidth]{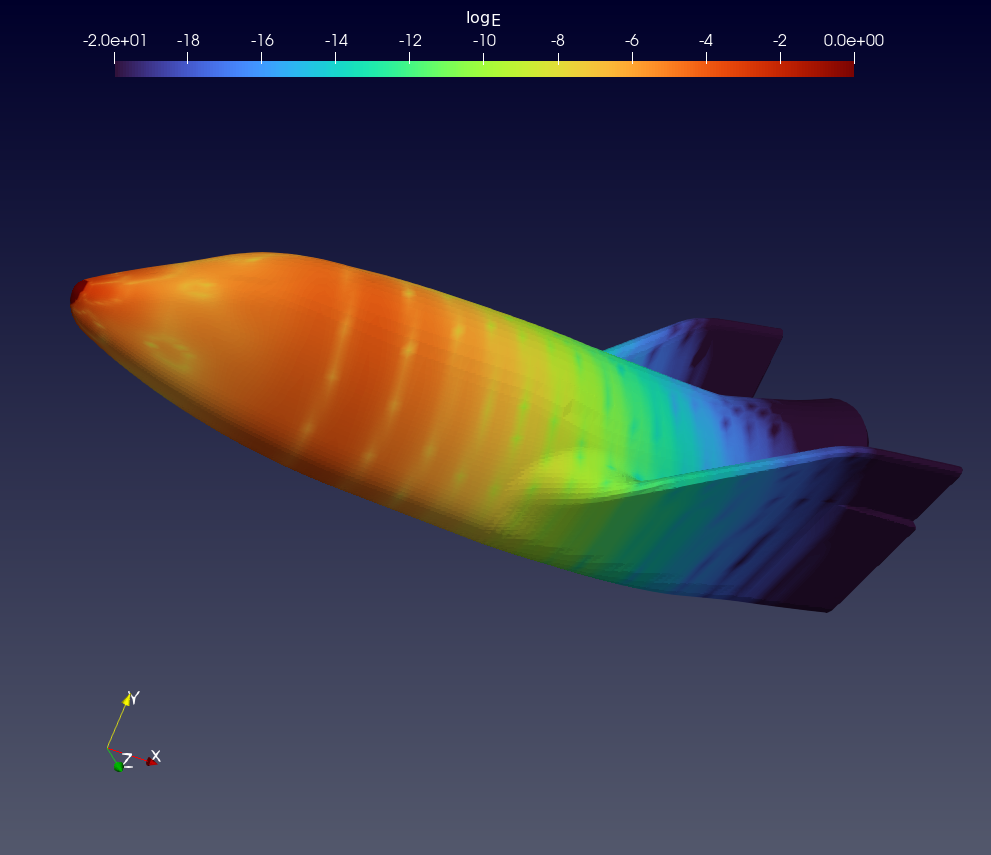}
    \caption{}
\end{subfigure}
\hfill
\begin{subfigure}[b]{0.48\textwidth}
    \includegraphics[width=0.95\textwidth]{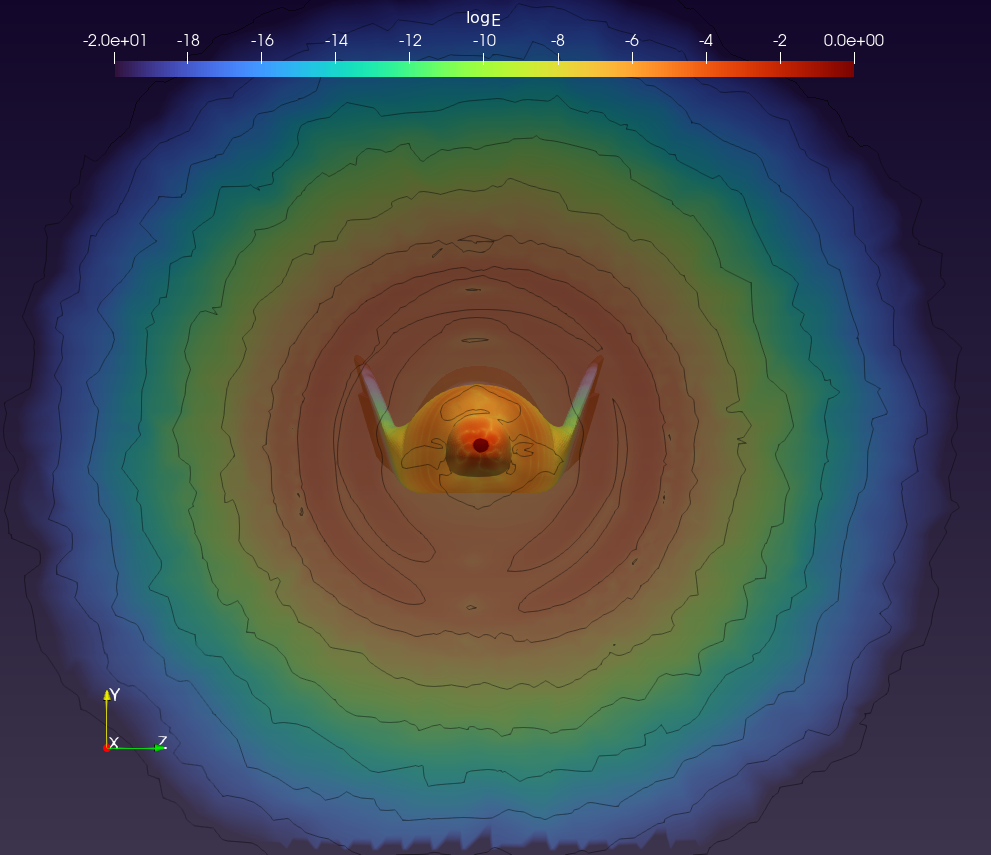}
    \caption{}
\end{subfigure}
\vfill
\begin{subfigure}[b]{0.48\textwidth}
    \includegraphics[width=0.95\textwidth]{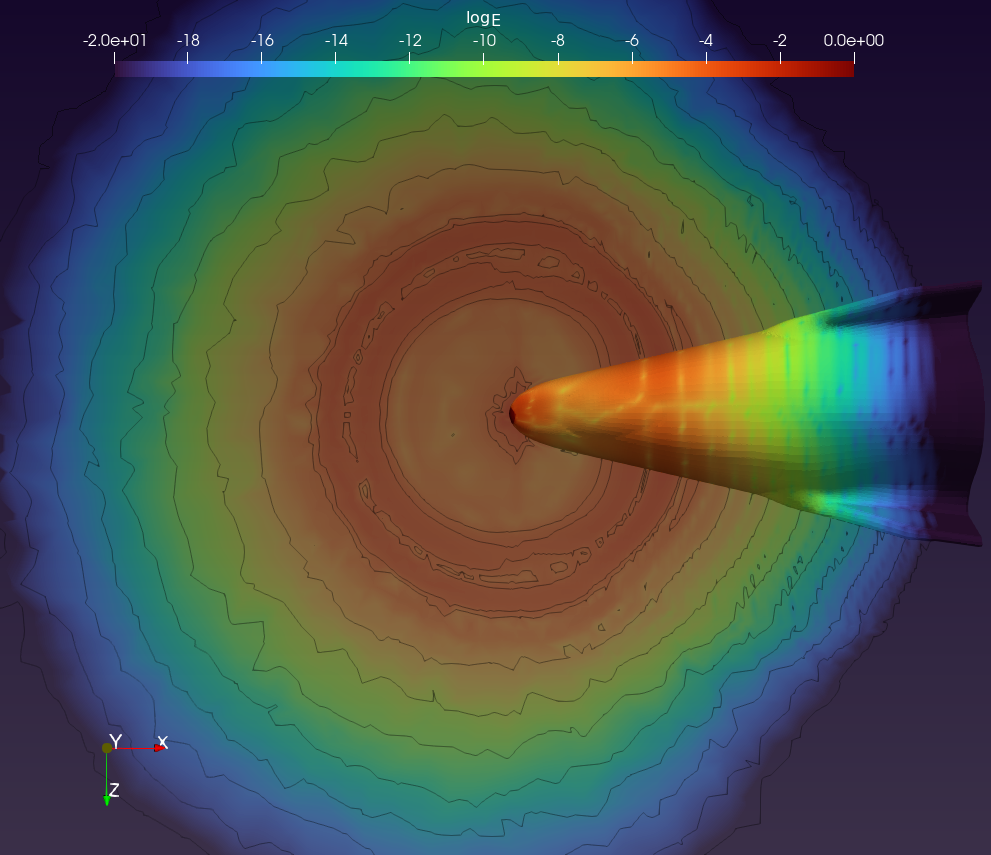}
    \caption{}
\end{subfigure}
\hfill
\begin{subfigure}[b]{0.48\textwidth}
    \includegraphics[width=0.95\textwidth]{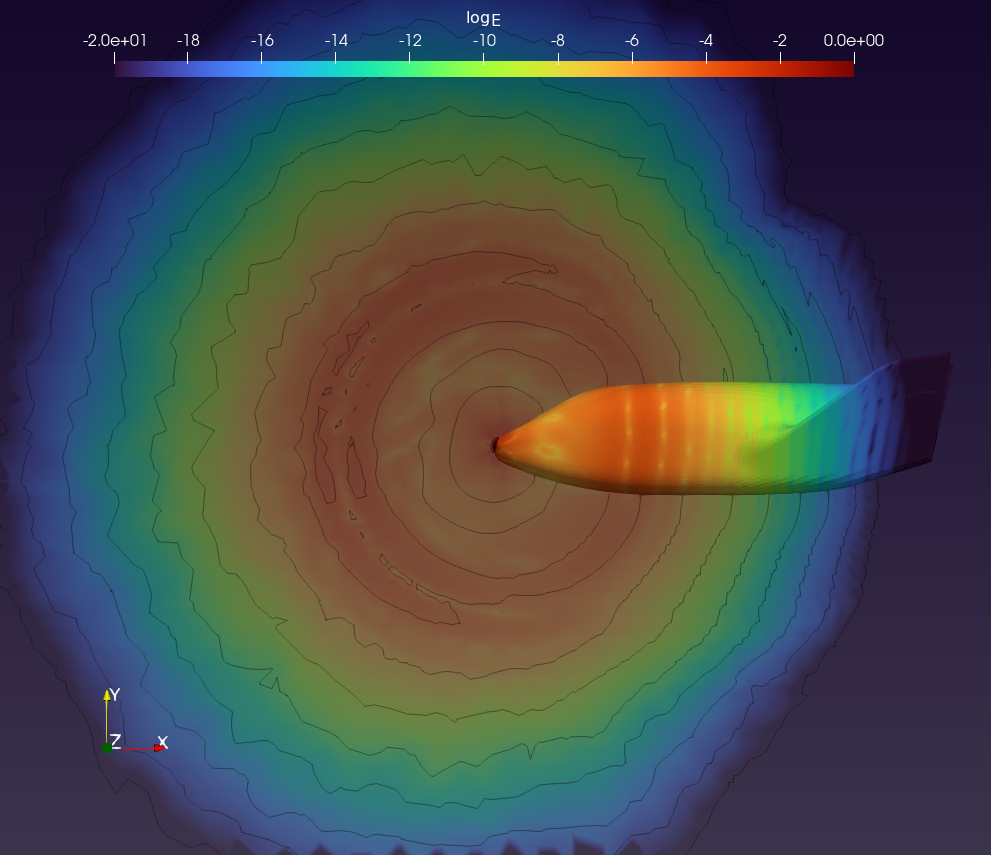}
    \caption{}
\end{subfigure}
\caption{Signal Propagation Around an Aircraft: The electric field magnitude on the surface of  (a) aircraft, (b) y-z slice, (c) x-z slice, (d) x-y slice.}
\label{fig:x382contour}
\end{figure}

\begin{figure}
\centering
\begin{subfigure}[b]{0.48\textwidth}
    \includegraphics[width=0.95\textwidth]{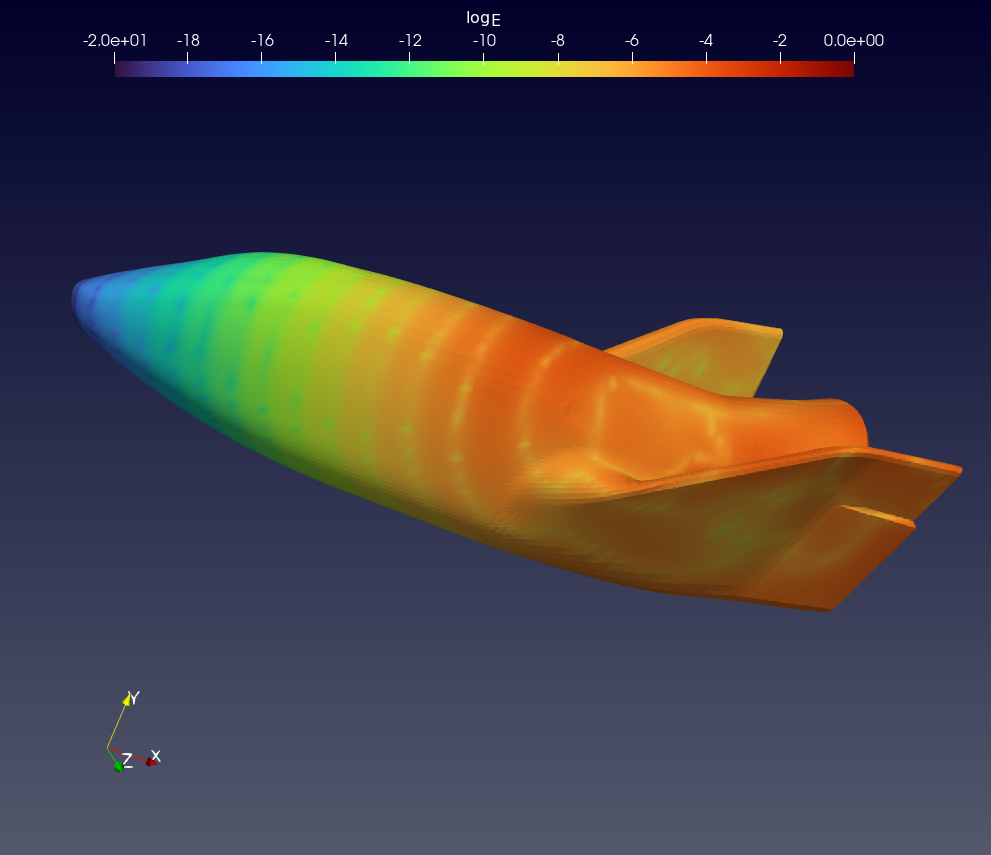}
    \caption{}
\end{subfigure}
\hfill
\begin{subfigure}[b]{0.48\textwidth}
    \includegraphics[width=0.95\textwidth]{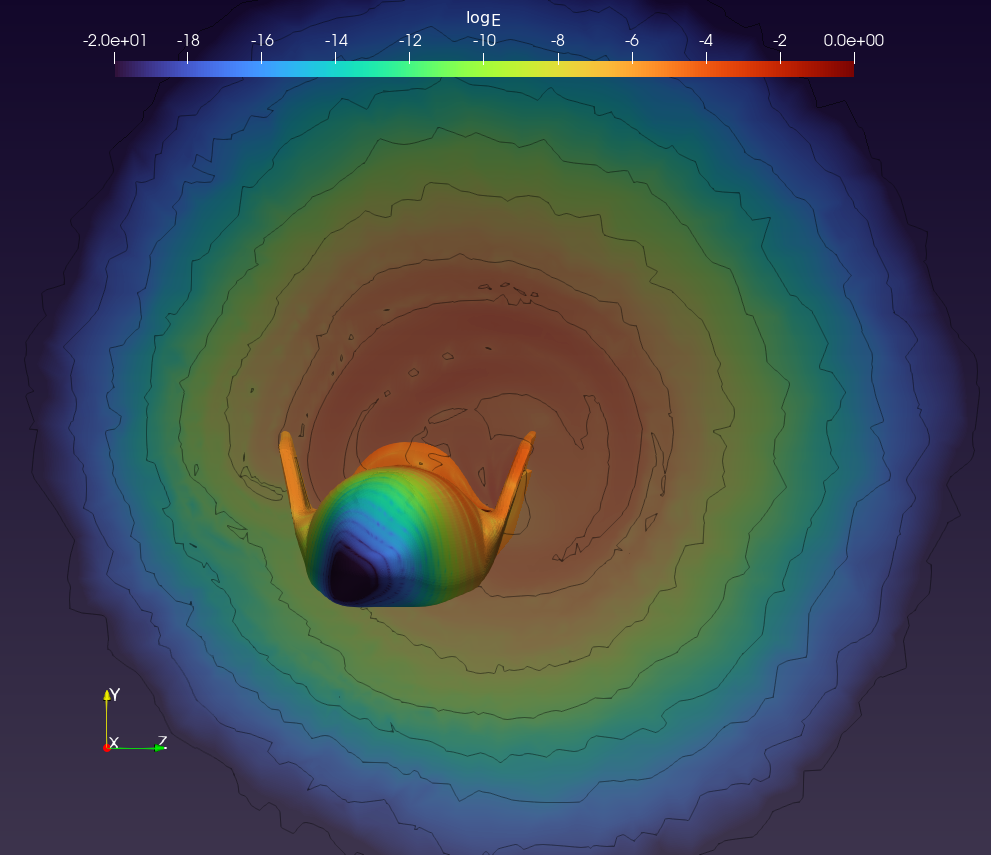}
    \caption{}
\end{subfigure}
\vfill
\begin{subfigure}[b]{0.48\textwidth}
    \includegraphics[width=0.95\textwidth]{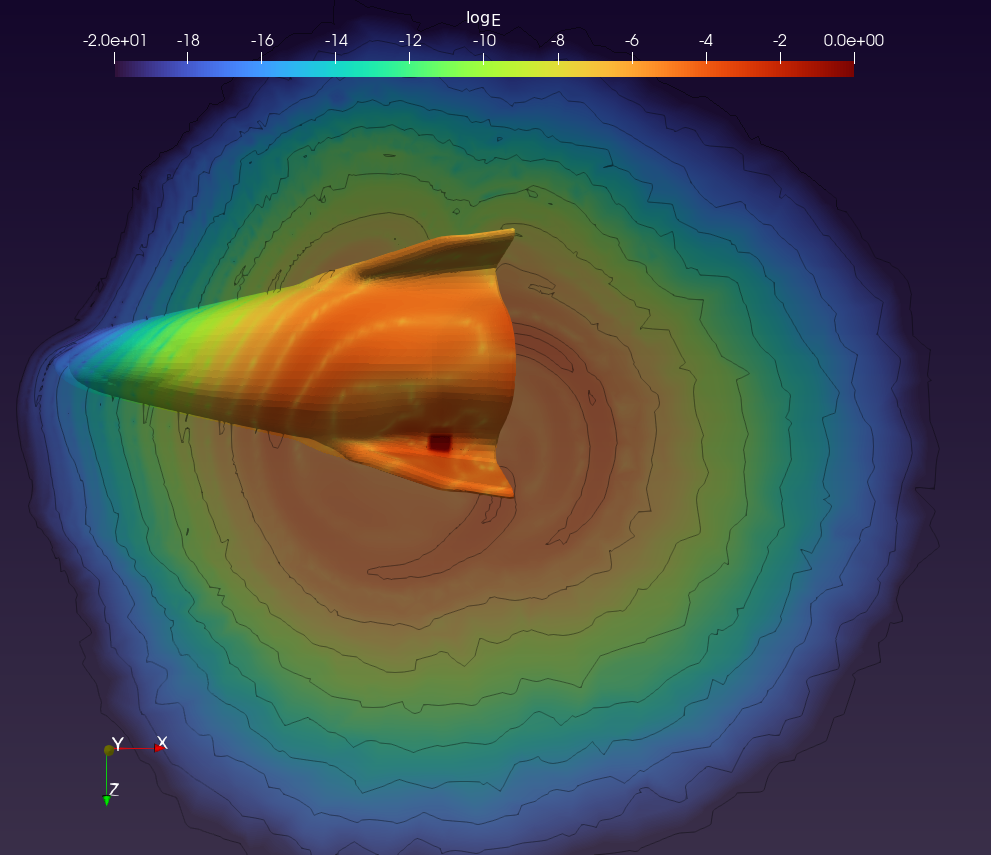}
    \caption{}
\end{subfigure}
\hfill
\begin{subfigure}[b]{0.48\textwidth}
    \includegraphics[width=0.95\textwidth]{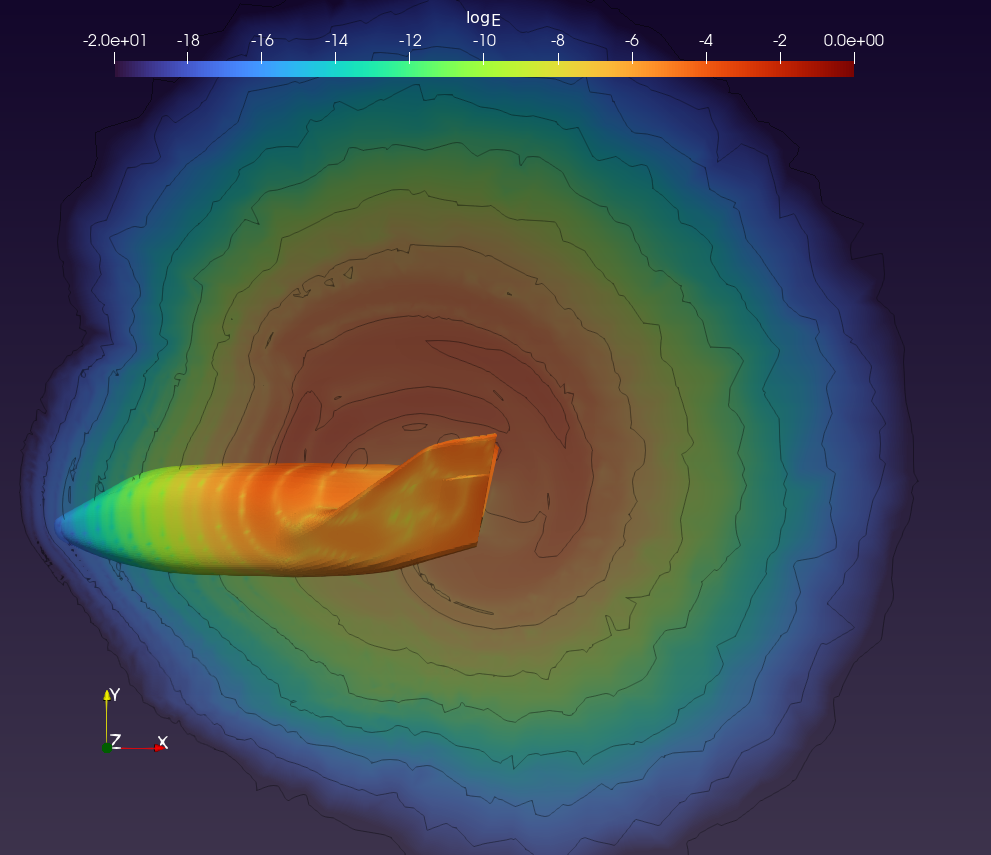}
    \caption{}
\end{subfigure}
\caption{Signal Propagation Around an Aircraft: The electric field magnitude on the surface of  (a) aircraft, (b) y-z slice, (c) x-z slice, (d) x-y slice.}
\label{fig:x383contour}
\end{figure}

Figures \ref{fig:x381contour}, \ref{fig:x382contour}, and \ref{fig:x383contour} present the simulation results. To visualize the rapid three-dimensional decay of the electric field, we plot the base-10 logarithm of its magnitude, \(\log_{10} ||\boldsymbol{E}||\). These figures illustrate the field distribution both on the aircraft surface and in the surrounding volume.

The specific antenna locations for each case are as follows:
\begin{itemize}
    \item Figure \ref{fig:x381contour}: The antenna is positioned at the center of the upper fuselage at coordinates (150, 73, 10) m. Its location is indicated by a localized region of very high field intensity.
    \item Figure \ref{fig:x382contour}: The antenna is mounted on top of the nose section at (5, 33, 2) m.
    \item Figure \ref{fig:x383contour}: The antenna is located within the fairing between the fuselage and empennage at (260, 42, 52) m.
\end{itemize}

The radiation patterns reveal distinct propagation characteristics influenced by the aircraft geometry. In Figure \ref{fig:x381contour}(b), signal propagation is more extensive in the upward direction than downward, owing to blockage by the aircraft body beneath the antenna. Figure \ref{fig:x381contour}(c) shows that propagation is also impeded in directions obstructed by the empennage. A similar trend is observed in Figure \ref{fig:x381contour}(d), where signal propagation is weakest in the lower-right quadrant.

In contrast, Figure \ref{fig:x382contour} demonstrates a more symmetric and homogeneous radiation pattern. With the antenna placed on the nose, the signal is less obstructed by the main airframe. Conversely, the pattern in Figure \ref{fig:x383contour} is highly asymmetric and inhomogeneous. This phenomenon results from the antenna's offset within the fairing and the significant scattering by nearby structural components.

Collectively, these results demonstrate that the GKS effectively captures the complex phenomena of electromagnetic wave propagation and scattering from a metallic aircraft surface.

\section{Conclusion}
\label{sec:conclusion}

This work has successfully extended GKS to computational electromagnetics using discrete velocity space. The proposed scheme achieves second-order accuracy in a single step by a time-accurate numerical flux at cell interfaces. Its inherent kinetic foundation provides a multidimensional framework, while the finite-volume discretization ensures natural extension to unstructured meshes. The inherent collision process gives lower numerical dissipation than classical FVS methods, and the formulation permits the straightforward implementation of non-reflecting boundary conditions.
The scheme is rigorously validated against standard electromagnetic benchmarks and compared with established methods, including FDTD, FVS, and LBM. Finally, its capability for modeling complex, realistic geometries is demonstrated through a large-scale simulation of electromagnetic wave propagation around a complete aircraft configuration. This extension of GKS establishes a robust, accurate, and flexible framework for electromagnetic simulations, particularly in applications involving complex geometries.

\section*{Acknowledgements}

This work was supported by the National Key R\&D Program of China (Grant No. 2022YFA1004500) and the National Natural Science Foundation of China (Nos. 12172316 and 92371107), and the Hong Kong Research Grant Council (Nos. 16208021, 16301222, and 16208324).\\

\textbf{Declaration of interests}: The authors report no conflict of interest.

\bibliographystyle{elsarticle-num}
\bibliography{ref.bib}

\end{document}